\documentclass[12pt]{article}

\usepackage[normalem]{ulem}
\usepackage{fancyhdr}
\usepackage{amsmath,amsfonts,amssymb}
\usepackage{graphicx}
\usepackage{color}
\usepackage[footnotesize]{caption}
\usepackage{epstopdf}
\usepackage[Lenny]{fncychap}
\usepackage{algorithm}
\usepackage{algpseudocode}
\usepackage{tocvsec2}
\usepackage{mathtools}
\usepackage{nicefrac}
\usepackage{tikz}
\mathtoolsset{showonlyrefs}

\makeatletter

\newcommand{\Rmnum}[1]{\expandafter\@slowromancap\romannumeral #1@}
\makeatother

\newtheorem{remark}{Remark}

\def\methone{ \textsc{Method 1} }
\def\methtwo{ \textsc{Method 2} }
\def\meththree{ \textsc{Method 3} }

\def\R{{\mathbb R}}

\let\theta\vartheta

\def\vec#1{{\mathbf{#1}}}

\DeclareMathOperator{\sgn}{sign}

\title{High order finite-difference ghost-point methods\\
for elliptic problems in domains with curved boundaries}
\author{Armando Coco, Giovanni Russo}

\begin{document}

\maketitle


\abstract{In this paper a fourth order finite difference ghost point method for the Poisson equation on regular Cartesian mesh is presented. The method can be considered the high order extension of the second ghost method introduced earlier by the authors. Three different discretizations are considered, which differ in the stencil that discretizes the Laplacian and the source term.
It is shown that only two of them provide a stable method. 
The accuracy of such stable methods are numerically verified on several test problems. 
}

\section{Introduction}
Elliptic equations in two space dimensions are often encountered in several physical and mathematical contexts, ranging from electrostatic (the Poisson equation describes the  electric potential associated to a distribution of charges) to electromagnetic (the Helmoltz equation models electromagnetic wave propagation in the frequency domain), from incompressible gas dynamics (the pressure has to satisfy an elliptic equation to impose that the flow is incompressible) to linear elasticity. Because of its ubiquitous applications, solution of elliptic problems has always been of great interest, from modeling, analytical and numerical points of view. 

The most classical approach to the solution of such equations is usually based on finite difference or finite element methods. 
Typically, finite difference schemes are extremely simple to implement producing efficient codes when discretizing the equations on a regular Cartesian grid on a rectangular domain. For non rectangular domains, the finite element approach has the advantage of a greater flexibility, since it can more easily treat complex geometries by adopting, for example, a triangulation which well approximates the boundary of the domain. 

More recently, advancements have been made in the development of finite difference methods on regular Cartesian grids. These methods discretize elliptic equations within arbitrary domains implicitly defined, for instance, by a level set function~\cite{CocoRusso:Elliptic, coco2018second, gibou2002second, Gibou:guidelinesPoisson, Osher-Fedkiw:level_set, Zhou:MIB, Iollo:penalization}.

Such schemes have the advantage offered by regular Cartesian mesh, and are very efficient for the treatment of problems with time dependent boundaries, since the grid does not need to be re-computed at every time step. Shortley-Weller method~\cite{Shortley-Weller:discretization} for the Poisson equation is based on the use of additional grid points at the intersection of the boundary with grid axis, and on the use of such additional points for the discretization of the Laplacian operator near the boundary (see left panel in Fig.~\ref{fig:classic}). It can be shown (see~\cite{Shortley-Weller:discretization}) that even if the discrete operator obtained in this way is only first order accurate, the resulting scheme based on the five-point approximation of the Laplacian applied to a Dirichlet problem guarantees second order accuracy for the solution.

A different alternative is offered by {\em ghost point methods}, which are  based on the following idea. 
The computational domain $\Omega$ is immersed in a rectangular region $R$, which is discretized by regular Cartesian mesh. The grid points on the mesh are classified in three disjoint sets: internal points (points inside $\Omega$), ghost points (points outside $\Omega$ with a neighbor inside $\Omega$) and inactive points (see right panel of Figure \ref{fig:classic}).
The unknown function is defined only on internal and ghost points. The discrete equations are defined on the internal points, while boundary conditions are assigned by imposing suitable conditions on the ghost points. 
In \cite{CocoRusso:Elliptic}, for each ghost point $G$ the boundary condition is assigned on the closest point $B_G$ on the boundary, while in~\cite{gibou2002second} the boundary conditions are assigned at the intersection between the boundary of $\Omega$ and the grid axis, giving rise to two different values of the function on the ghost point.

More recently, a different version of ghost based method defined on a regular Cartesian mesh for elliptic problems has been proposed \cite{astuto2024nodal}, which is second order accurate and guarantees better theoretical properties (symmetry, coercivity) leading to a convergence proof under mild regularity assumptions on the boundary of the domain $\Omega$. 

In this paper we propose a fourth order extension of the second order method introduced in \cite{CocoRusso:Elliptic} for the Poisson equation in two dimensions. 
Other high order finite-difference methods for Elliptic problems are~\cite{fernandez2020very, baeza2016high, ren2022fft, gabbard2024high, feng2022high}.

High order methods  are more cost-effective when accurate solutions are needed. The benefit appears when the solution itself is sufficiently regular, which implies regularity of the data (i.e.\ regularity of the boundary $\partial \Omega$, the data on the boundary and the right hand side of the equation). For free boundary problems, in which the boundary of the domain $\partial \Omega$ evolves according to an assigned level set equation, there has been an increase interest in fourth order schemes for the evolution of the level set function and for the computation of the corresponding signed distance function \cite{Gibou:reinizialization}, which allows a second order accurate computation of the surface curvature. 

\subsection{Model problem}
\begin{figure}
\begin{center}
\includegraphics[width=0.5\textwidth]{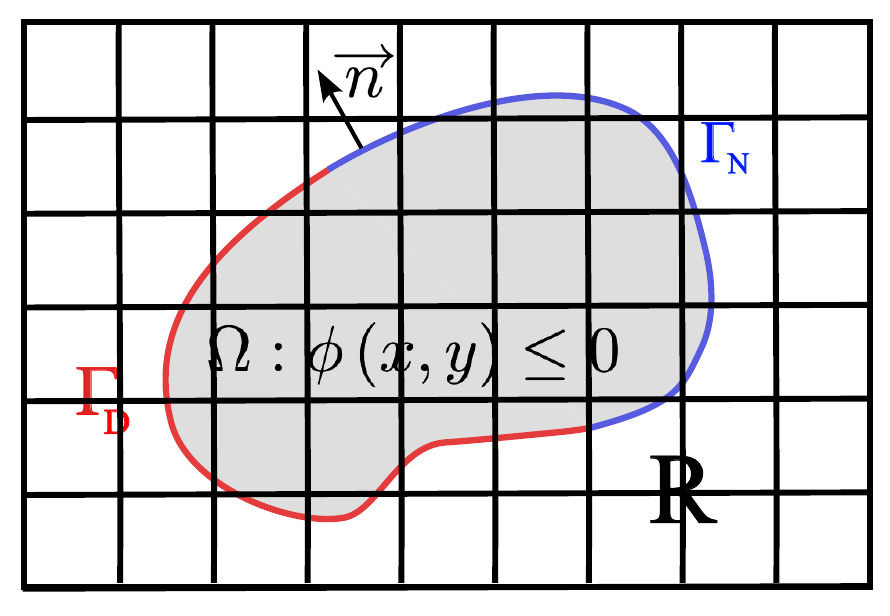}
\end{center}
\caption{Setup of problem \eqref{eq:modprob}: the domain $\Omega$, defined by $\phi(x,y)\leq 0$, is embedded into a rectangular domain $R$. Mixed boundary conditions are imposed on the frontier $\Gamma$ of the domain. In each point of $\Gamma$ the unit normal pointing outside of the domain can be expressed as $\vec{n}=\nabla\phi/|\nabla\phi|$. }
\label{fig:Omega}
\end{figure}

We consider the following mixed Dirichlet-Neumann problem for for the Poisson equation
\begin{equation}\label{eq:modprob}
\begin{split}
- \Delta u  &= f \text{ in } \Omega \\
u &= g_D \text{ on } \Gamma_D \\
u &= g_N \text{ on } \Gamma_N \\
\end{split}
\end{equation}
The compact domain $\Omega$ is a subset of a rectangle $R$ of $\mathbb{R}^2$, i.e.\ $\Omega\subset R\subset \mathbb{R}^2$, and it is  implicitly defined by a level set function $\phi:R\to \mathbb{R}$, $\Omega = \{(x,y)\in R | \phi(x,y)\leq 0\}$. 
We assume that the boundary of the domain $\Gamma = \partial \Omega$, defined by 
$\Gamma = \{(x,y)\in R \, | \, \phi(x,y) = 0\}$, is sufficiently smooth, which is guaranteed if the function $\phi$ is sufficiently regular in a neighborhood of $\Gamma$, and its gradient is bounded away from zero. 
The boundary $\Gamma$ is partitioned in a Dirichlet boundary $\Gamma_D$ and  Neumann boundary $\Gamma_N$: $\Gamma = \Gamma_D\cup \Gamma_N$, $\Gamma_D\cap\Gamma_N = \emptyset$.
As usual, we denote by 
$\Delta \equiv \partial^2/\partial x^2 + \partial^2/\partial y^2$
the Laplacian operator in $\mathbb{R}^2$.
The method proposed in this paper is high order accurate for sufficiently regular source term and boundary data, therefore we assume that $f$, $g_D$ and $g_N$ are $\mathcal{C}^1$ functions in $\Omega$.
The problem setup is illustrated in Fig.~\ref{fig:Omega}.

\subsection{Review of the second order method}
\begin{figure}
\begin{center}
\begin{minipage}{0.40\textwidth}
\includegraphics[width=0.9\textwidth]{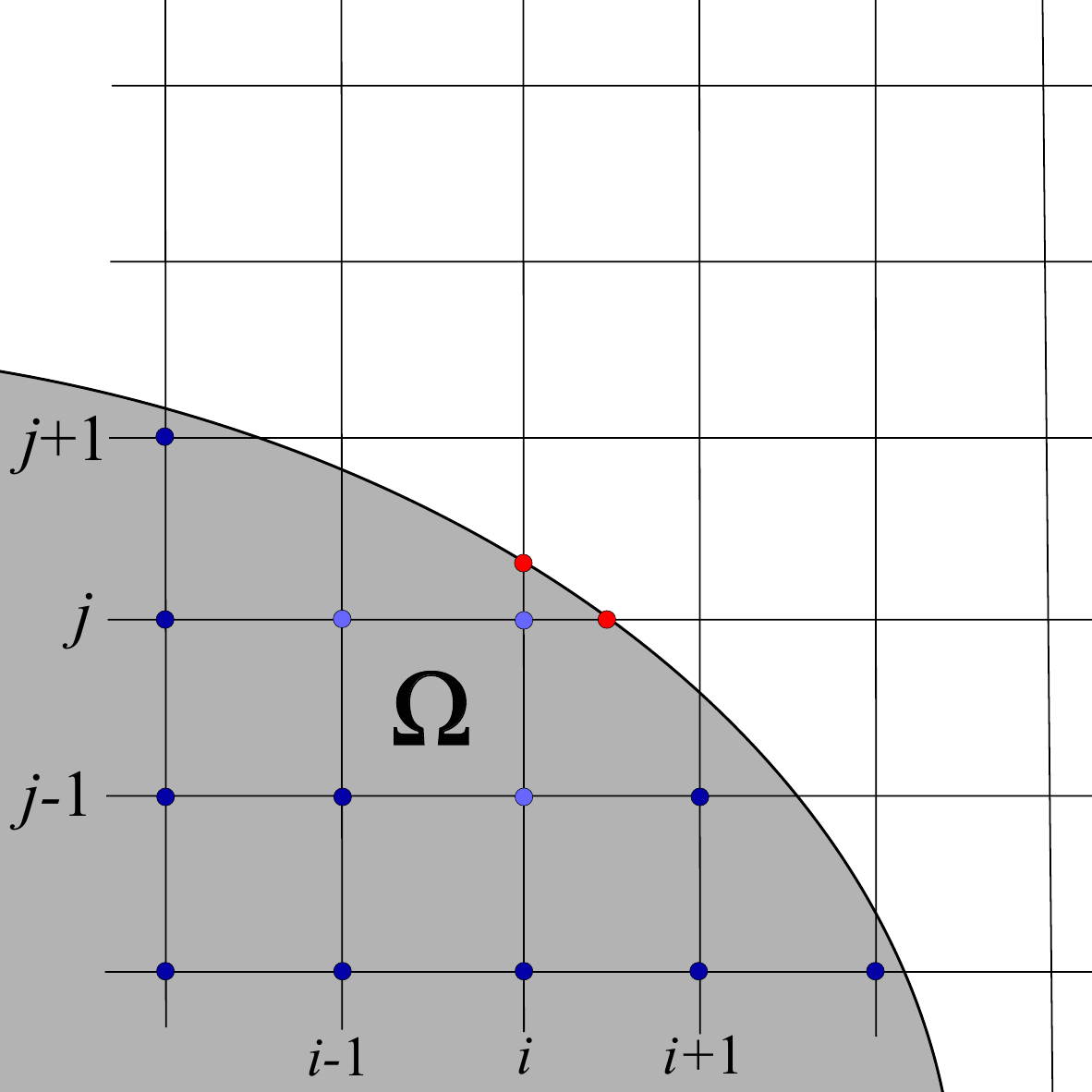}
\end{minipage}
\begin{minipage}{0.40\textwidth}
\includegraphics[width=0.9\textwidth]{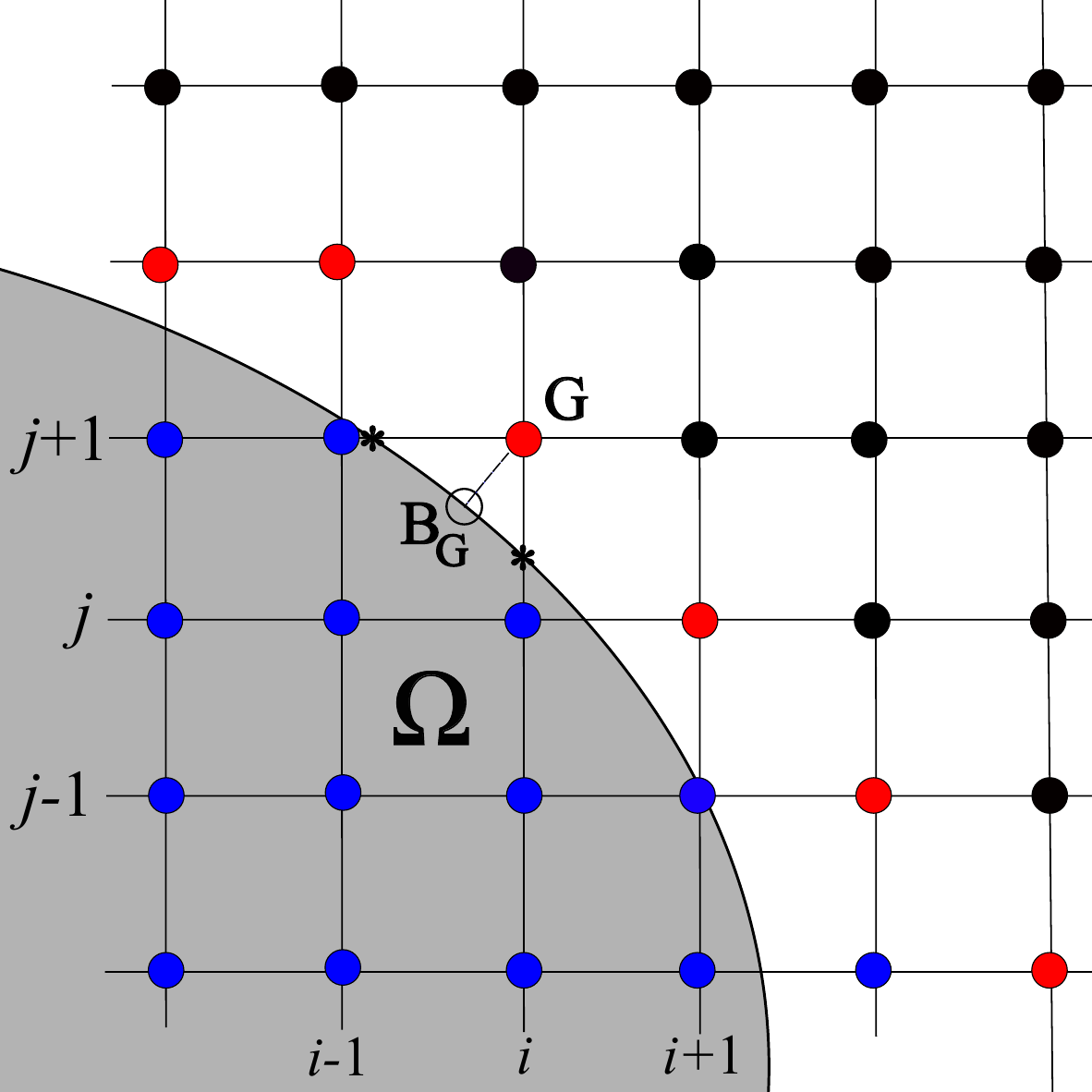}
\end{minipage}
\end{center}
\caption{Left panel: Shortley-Weller stencil around grid point $(i,j)$. Right panel: Ghost point discretization for the classic five-point discrete Laplacian. Boundary condition on node $G$ are assigned by imposing that polynomial that interpolated the solution on the (internal of ghost) nodes sarisfied the boundary conditions on point $B_G$, as in \cite{CocoRusso:Elliptic}, or on intersection of the grid lines with the boundary (denoted by the asteriscs), as in~\cite{Gibou:Ghost}.}
\label{fig:classic}
\end{figure}

Here we give a brief review of the second order methods presented in \cite{CocoRusso:Elliptic}.
The rectangular region is meshed by a regular Cartesian square grid, and throughout the paper we assume that the grid spacing $h=\Delta x = \Delta y$ is considerably smaller than the smallest radius of curvature of the boundary $\Gamma$, although the method can be extended to the case of Lipschitzian boundaries by adopting the approach presented, for example, in~\cite{CocoRusso:Elliptic, Zhou:MIB}.
In the internal points, the Laplacian is discretized by the classical second order 5-point discrete Laplacian. 
For each ghost point, a condition is assigned by imposing that a suitable reconstruction of the solution satisfies the boundary condition on the projection of the ghost point on the boundary (see right panel of Fig.~\ref{fig:classic}. The simplest scheme is obtained by using bilinear 
interpolation\footnote{Here by bilinear, biquadratic, bicubic, etc.\ interpolation we mean interpolation obtained on a regular square grid, using multivariate polynomials in two space variables, of degree at most $n$, with $n$ respectively, equal to 1, 2, 3, etc.\ in each variable. Such multivariate polynomials will take the form $p_n(x,y) = \sum_{i=0}^n\sum_{j=0}^n a_{ij}x^iy^j$. The name comes from the consideration that each such interpolating polynomial can be obtained by repeated interpolations in the coordinate directions by one dimensional polynomials of degree at most $n$.} For example, with reference to the right panel of Fig.~\ref{fig:classic}, assuming we have Dirichlet boundary conditions near ghost point $G$, the condition is obtained by imposing that the bilinear reconstructions from the values in the four points indexed by $(i,j),(i,j+1),(i-1,j+1),(i-1,j)$, evaluated in $B_G$, is equal to $g_D(B_G)$. Notice that in the particular case in the figure, the condition in point $G$ is independent on the conditions on the other ghost points, so the value imposed in $G$ could be expressed in terms of the internal points. However, in general this is not the case (as, for example, for point $(i-1,j+2)$), so at the end one obtains a coupled system for the internal and ghost values. 
It has been shown that this simple scheme provides second order accuracy for the function $u$ in the case of Dirichlet boundary conditions, and first order accuracy for Neumann conditions. Once the function $u$ is known on internal and ghost points, the gradient can be easily computed on the internal nodes by standard centered discretization. However the technique based on bilinear interpolation only provides first order accuracy for the evaluation of the gradient.
A better accuracy is obtained by using biquadratic interpolation, which is based on a 9-point stencil pointing in the direction opposite to the normal, direction that we denote as {\em upwind\/}. In the case of point $G$ in the figure, for example, the biquadratic interpolation uses points $(i-\ell,j+1-k)$, $\ell=0,1,2$ and $k=0,1,2$.
Such stencil guarantees second order accuracy both for the function $u$ and its gradient.
In order to improve the efficiency in the solution of the resulting linear system for the internal and ghost values, a geometric multigrid method has been proposed, which is based on a suitable smoother constructed {\it ad hoc\/} for the ghost point discretization. 

The method has been generalized to the case in which the boundary $\Gamma$ is only piecewise regular~\cite{CocoRusso:Elliptic} and to the case of a more general elliptic equation with piecewise continuous coefficients \cite{coco2018second}. The method has been successfully applied to problems in linear elasticity, with application to volcanology \cite{coco2014second}. 
The case of moving boundaries has been considered in \cite{coco2020multigrid}.
An improved multigrid applied to the ghost method has been proposed and analyzed in \cite{coco2023ghost}, by optimizing the relaxation parameters associated with the boundary condition to augment the method's smoothing efficacy.
Finally, such boundary treatment has been adopted for hyperebolic problems, both in the case of fixed and moving boundaries 
\cite{chertock2018second}.

The convergence of the method has been studied in \cite{coco2023spectral}, in the particular case of rectangular domains whose boundary does not fit with the Cartesian lines. A general proof of convergence of the ghost method developed in \cite{CocoRusso:Elliptic} for a general 2D domain has never been obtained. 

A recent cut-FEM variant of the ghost method has been presented in~\cite{astuto2024nodal}. The method results in a symmetric, positive definite and coercive linear system, for which it is possible to prove convergence for a general domain. 

The rest of the paper is structured as follows. After this introduction, in Section 2 we present the various space discretizations of the equation. The third Section is devoted to the treatment of ghost points and boundary conditions. Section 4 illustrates numerical results on two test problems. In the final section we draw some conclusions. 

\section{Numerical Method}
The computational domain $\mathcal{D} = [-1,1]^2$ is discretized by a uniform Cartesian grid of $(n_x+1) (n_y+1)$ equally spaced grid points $(x_i,y_j)$, with $x_i=-1+ih_x$, $y_j=-1+jh_y$ for $i=0, \ldots, n_x$, $j=0, \ldots, n_y$, where $h_x = 2/n_x$ and $h_y=2/n_y$ are the spatial steps in the $x$ and $y$ directions, respectively. For simplicity, we assume that $n_x=n_y=N$, then $h_x=h_y=h$, but the entire methodology can be easiliy extended to the general case $h_x \neq h_y$. The set of grid points is denoted by $\mathcal{D}_h$. Grid points that are inside $\Omega$ are called {\em internal grid points}, forming the discrete version of $\Omega$, denoted by $\Omega_h = \mathcal{D}_h \cap \Omega$. The remaining points are called 
{\em external grid points}.
We use the notation $\omega_{i,j}$ to denote the approximation of $\omega(x_i,y_j)$ for any function $\omega \colon \Omega \rightarrow \mathbb{R}$.

\subsection{Discretization of the elliptic equation}
Several fourth order finite difference discretization of \eqref{eq:modprob} are possible~\cite{fernandez2020very, baeza2016high, ren2022fft, gabbard2024high, feng2022high}.
In this paper, we examine three distinct discretizations based on a fully coupled ghost-point approach, denoted as \methone, \methtwo and \meththree. We demonstrate that the first discretization is too ill-conditioned to ensure convergence. In contrast, the other two exhibit the expected fourth-order convergence rate.

The PDE (first equation of \eqref{eq:modprob}) is approximated in each internal grid point using a fourth-order discretization. We investigate two types of discretizations: star-stencil and box-stencil discretizations. 
We use the following matrix notation to denote the linear combination of grid values  $\omega_{i,j}$.

 \begin{equation}\label{eq:matrixnotationiC}
 S\omega_{i,j}
 =
 \sum_{p,q=-P}^P s_{p,q} \omega_{i+p,j+q}
 \end{equation}
where $S\in\R^{(2P+1)\times(2P+1)}$, $P\in\{0,1,2\}$.
With this notation, the discretization of Eq. \eqref{eq:modprob} takes the form 
\[
    S\,u_{i,j} = Mf_{i,j} + \mathcal{O}(h^4)
\]
We consider the following two fourth order discretizations.
\subsubsection*{Star-stencil discretization}\label{sect:star}
The star-stencil discretization consists of a natural  higher order extension of the classical second-order accurate central difference five-point star stencil discretization. Second derivatives are discretized along the Cartesian coordinates  by using central difference formulas obtained by Taylor expansions up to the desired order of accuracy. The elliptic equation $- \nabla^2 u (x_i,y_j) = f(x_i,y_j)$ is dscretized by:
\begin{equation}\label{disc:star}
\frac{1}{12 h^2}
\begin{bmatrix}
& & 1 & &  \\
& & -16 & &  \\
1& -16 & 60 & -16 & 1 \\
& & -16 & &  \\
& & 1 & &  \\
\end{bmatrix}
u_{i,j}
= f_{i,j} + \mathcal{O}(h^4)
\end{equation}
for each internal grid point $(x_i,y_j) \in \Omega_h$, so $P=2$ for $S$ and $M = 1$.
This scheme is obtained by approximating the second derivative in each coordinate direction by the corresponding fourth order five-point discretization. 

\subsubsection*{Box-stencil (Mehrstellen) discretization}\label{sect:box}
The Mehrstellen discretization~\cite{collatz2012numerical} uses a compact $3 \times 3$ stencil around the internal grid point $(x_i,y_j)$. The elliptic equation $- \nabla^2 u (x_i,y_j) = f(x_i,y_j)$ is discretized as:
\begin{equation}\label{disc:box}
\frac{1}{6 h^2}
\begin{bmatrix}
-1 & -4 & -1 \\
-4 & 20 & -4 \\
-1 & -4& -1 \\
\end{bmatrix}
u_{i,j}
=
\frac{1}{12}
\begin{bmatrix}
 & 1 &  \\
1 & 8 & 1 \\
 & 1&  \\
\end{bmatrix} 
f_{i,j} + \mathcal{O}(h^4)
\end{equation}
for each internal grid point $(x_i,y_j) \in \Omega_h$, so $P=1$ both for $S$ and $M$.
This scheme can be understood as follows. It can be shown that matrix $S$ provides a fourth 
order approximation of the operator $-\Delta -h^2 \Delta^2$, i.e. 
\[
    S u(x_i,y_j) = -\Delta u(x_i,y_j) - h^2 \Delta^2 u(x_i,y_j) + O(h^4)
\]
while the matrix $M$ is a fourth order approximation of the operator $I+h^2\Delta$, i.e.
\[
    M f(x_i,y_j) = f(x_i,y_j)+h^2 Lf(x_i,y_j) = f(x_i,y_j) + h^2 \Delta f(x_i,y_j) + O(h^4)
\]
where
\[
    L = \frac{1}{h^2}
    \begin{bmatrix}
      & 1 &  \\
     1 & -4 & 1 \\
      & 1&  \\
     \end{bmatrix} 
\]
denotes the matrix corresponding to the classic five point stencil Laplacian.
The approximation is then obtained by observing that 
\[
    - \Delta u - h^2 \Delta^2 u = f + h^2\Delta (\Delta u) = f + h^2 \Delta f
\]

\subsection{Domain representation by level-set}
The physical domain $\Omega$ is implicitly described by a level-set function $\phi : \mathcal{D} \rightarrow \mathbb{R}$ such that:
\[
\Omega = \left\{ (x,y) | \phi(x,y) < 0 \right\}
\]
\[
\Gamma = \partial \Omega = \left\{ (x,y) | \phi(x,y) = 0 \right\}.
\]
Level set methods are largely used to track interfaces and they have been widely discussed in literature (see, for example, \cite{Osher-Fedkiw:level_set, Sethian:level_set}).
Different level-set functions can describe the same domain $\Omega$. For example, the disk centered at $(x_C,y_C)$ with radius $r>0$ is described by both the level-set functions $\phi_1 (x,y) = (x-x_C)^2 + (y-y_C)^2 - r^2$ and
$\phi_2 (x,y) = \sqrt{ (x-x_C)^2 + (y-y_C)^2} - r$. The latter is the so-called signed-distance function:
\[
\phi(x,y)=
\left\{
\begin{array}{rl}
 -d\left((x,y), \partial \Omega \right) & \mbox{ if } (x,y) \in \Omega, \\
 d\left((x,y),\partial \Omega\right) & \mbox{ if } (x,y) \notin \Omega,
\end{array}
\right.
\]
where $$d\left((x,y),\partial \Omega \right) =  \inf_{(\bar{x},\bar{y}) \in \partial \Omega} d_e\left( (x,y), (\bar{x},\bar{y}) \right)$$ is the distance between a point and the set $\partial \Omega$ ($d_e$ is the Euclidean distance between points).
For an arbitrary domain $\Omega$, the signed-distance function can be obtained from a generic level-set function $\phi_0(x,y)$, for example, 
by the reinitialization approach~\cite{SSO:level_set, Russo-Smereka:reconstruction, Gibou:reinizialization}, namely as the steady-state solution of the following evolutionary equation
\begin{equation}
\frac{\partial \phi}{\partial t} = \sgn(\phi_0) \left(1- \left| \nabla \phi \right| \right), \quad \phi=\phi_0 \; \mbox{ at } \; t=0,\label{eq:sdf}
\end{equation}
where $t$ is a fictitious time that represents an iterative parameter.
This equation updates the signed distance function away from the zero level set $\partial\Omega$, at a speed of one, conforming to the standard CFL-type stability conditions for the time step. Consequently, in a single time step, the signed distance function is updated by approximately one grid spacing.
The signed distance function is usually numerically more stable than a generic level-set function,
since the latter may develop steeper (or shallower) gradients, while $| \nabla \phi| =1$ for the signed-distance function.
We note that it is sufficient to solve Eq.~\eqref{eq:sdf} only for a few time steps, since we need to know the distance function only in the vicinity of the boundary.

The methods presented in this paper do not rely on the assumption that the level-set is a signed distance function, therefore they can be applied with generic level-set functions, provided the function is sufficiently smooth in a neighborhood of $\partial \Omega$.

From the level-set function we can recover geometric properties at each boundary point, such as the outward unit normal vector $\vec{n}$ (pointing from inside $\Omega$ towards outside) and the curvature $\kappa$:
\[
 \vec{n} = \frac{\nabla \phi}{\left| \nabla \phi \right|}, \; \; \; \kappa = \nabla \cdot \vec{n} = \nabla \cdot \frac{\nabla \phi}{\left| \nabla \phi \right|}.
\]

\section{Ghost-point method for boundary conditions}
\begin{figure}
\begin{center}
\begin{minipage}{0.40\textwidth}
\includegraphics[width=0.99\textwidth]{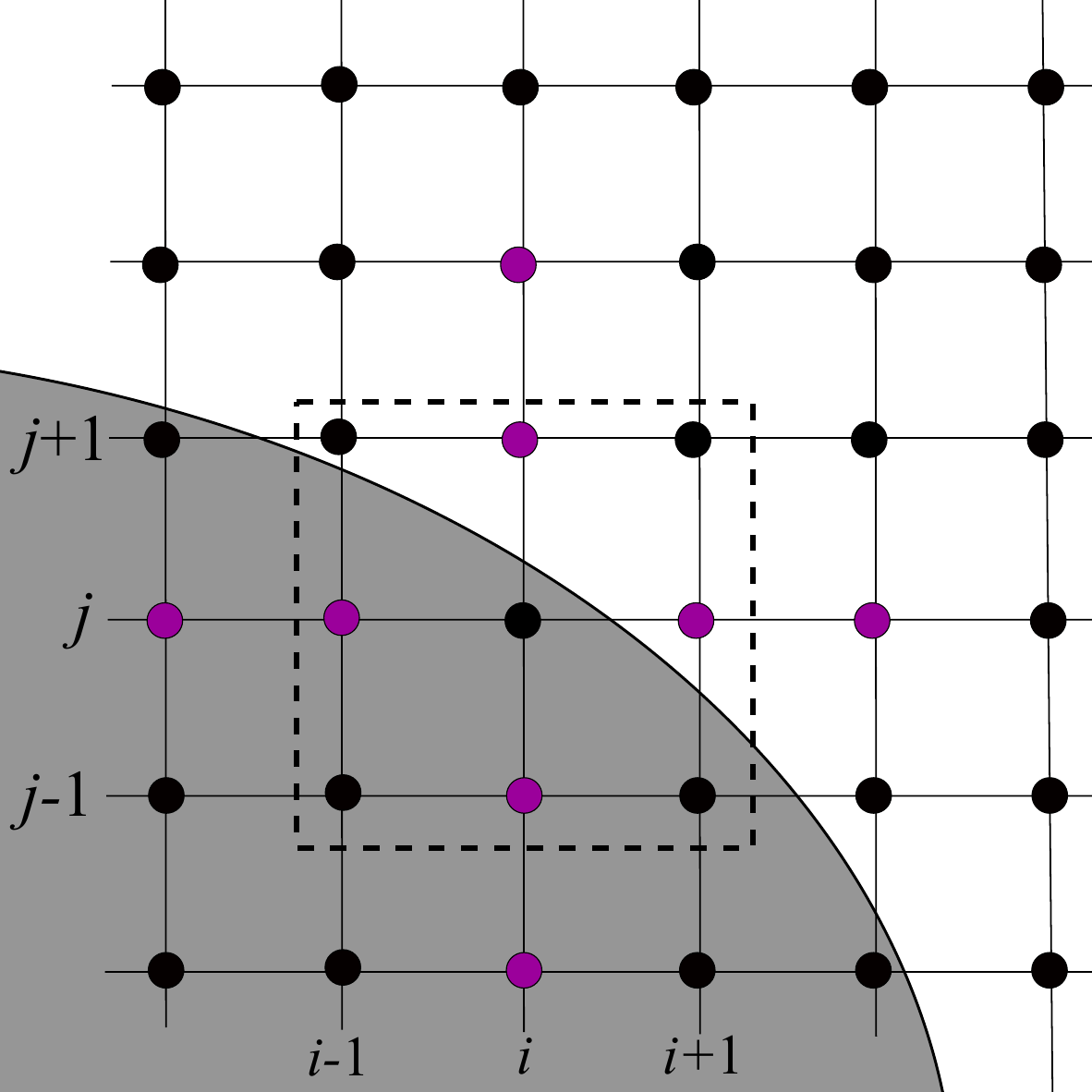}
\end{minipage}
\begin{minipage}{0.40\textwidth}
\includegraphics[width=0.99\textwidth]{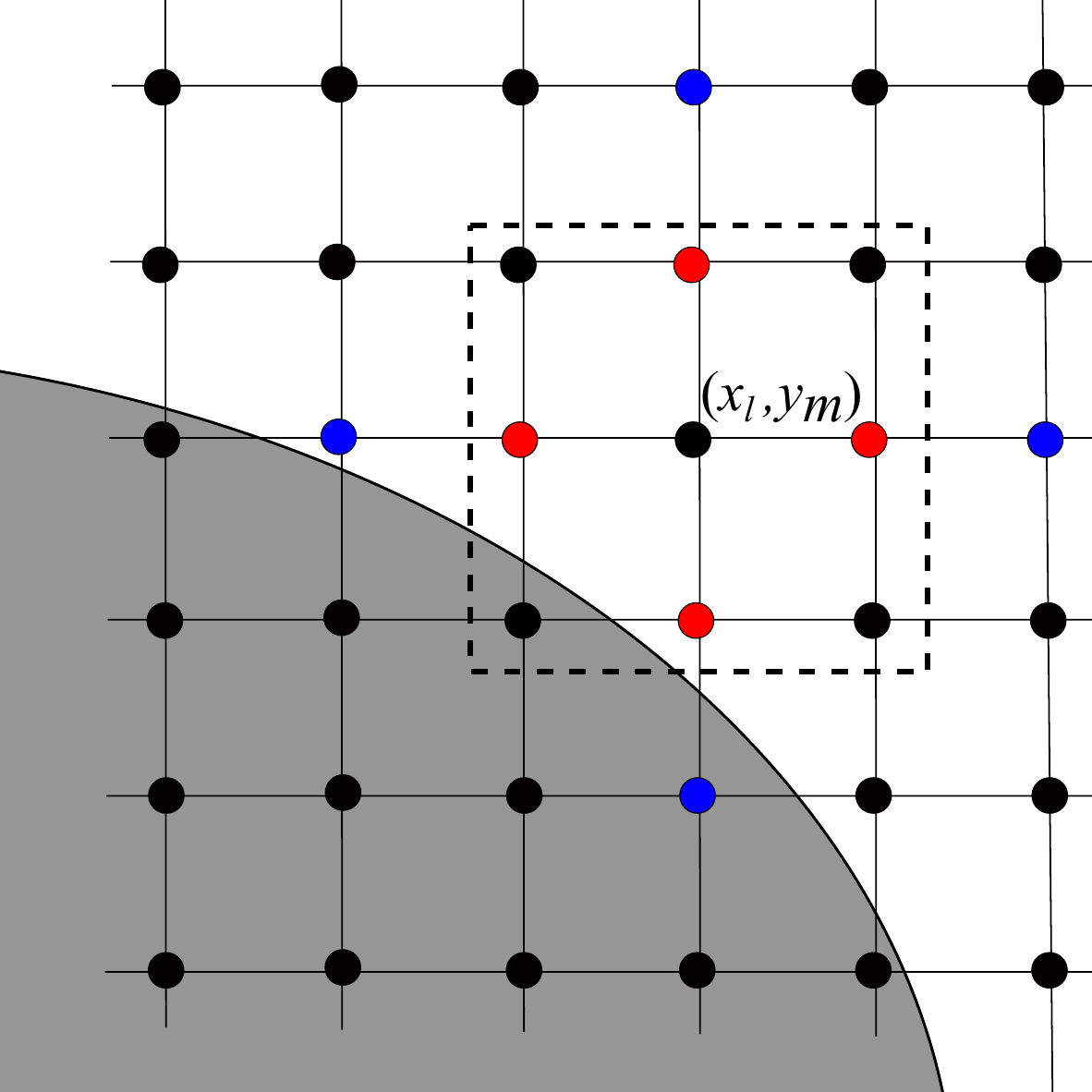}
\end{minipage}
\end{center}
\caption{
Left: Stencil adopted in the discretization of the Laplace operator on a generic internal grid point $(x_i,y_j)$. For the star-stencil discretization \eqref{disc:star}, the stencil is made by the purple grid points, while the dashed line surrounds the stencil of the box-stencil discretization \eqref{disc:box}.
Right: Star-shaped set of neighbour points $\mathcal{N}_{l,m}^{S,(1)}$ (red points) and $\mathcal{N}_{l,m}^{S,(2)}$ (blue points) for a generic ghost point $(x_l,y_m)$, as defined in Eq.~\eqref{eq:starNeigh}. The box-shaped set of neighbour points $\mathcal{N}_{l,m}^B$ (defined in Eq.~\eqref{eq:boxNeigh}) is surrounded by a dashed line. The set of neighbors of external grid points is adopted to identify the ghost points.}
\label{fig:stencil}
\end{figure}

When the internal grid point $(x_i,y_j)$ is sufficiently close to the boundary, the internal discretizations \eqref{disc:star} and \eqref{disc:box} require values of $u$ on external grid points, since the stencil is intersected with $\partial \Omega$. For example, in Fig.~\ref{fig:stencil} (left panel), the star-stencil discretization of the Laplace operator on $(x_i,y_j)$ requires the value of $u$ on purple grid points, including the external grid points
$(x_{i+1},y_j), (x_{i+2},y_j), (x_{i},y_{j+1}), (x_{i},y_{j+2})$, while the box-stencil discretization requires the value of $u$ on the points surrounded by a boxed dashed line, including the external points
$(x_{i-1},y_{j+1}), (x_{i},y_{j+1}), (x_{i+1},y_{j+1}), (x_{i+1},y_{j})$. The aim is to define a fictitious value of $u$ on external grid points ({\em ghost value}) by enforcing the boundary conditions in such a way that the overall fourth-order accuracy is not degraded. 

The external grid points that are required by the discretization of the Laplace operator are called {\em ghost points}. Depending on the type of the internal discretization, namely star- or box- stencil, we have two different sets of ghost points $\Gamma_h^{S}$ and $\Gamma_h^{B}$, respectively, defined as follows.
We first note that if a grid point $P$ is included in the discretization stencil of the Laplacian operator centered at another grid point $Q$, then $Q$ also appears in the stencil centered at $P$. Therefore, it is possible to define the set of ghost points by considering the stencil centered at the ghost point itself, rather than at an internal point. In detail,
for each external grid point $(x_l,y_m)$, we define the star-shaped set of neighbor points with grid distance $k\in\{1,2\}$:
\begin{equation}\label{eq:starNeigh}
\mathcal{N}_{l,m}^{S,(k)} = \left\{ (x_{l \pm k},y_m), (x_{l},y_{m \pm k}),  \right\} \cap \mathcal{D}_h
\end{equation}
and the box-shaped set of neighbour points:
\begin{equation}\label{eq:boxNeigh}
\mathcal{N}_{l,m}^B= \left\{ (x_{l + k_x},y_{m + k_y}), \; k_x,k_y =-1,0,1  \right\} \cap \mathcal{D}_h.
\end{equation}
The grid point $(x_l,y_m)$ is ghost if at least one of its neighbors is internal. Formally, the sets of ghost points are defined as:
\[
 (x_l,y_m) \in \Gamma_h^{S,(1)} \Longleftrightarrow \mathcal{N}_{l,m}^{S,(1)} \cap \Omega_h \neq \emptyset,
\]
\[
 (x_l,y_m) \in \Gamma_h^{S,(2)} \Longleftrightarrow (x_l,y_m) \notin \Gamma_h^{S,(1)} \text{ and } \mathcal{N}_{l,m}^{S,(2)} \cap \Omega_h \neq \emptyset,
\]
\[
 (x_l,y_m) \in \Gamma_h^{B} \Longleftrightarrow \mathcal{N}_{l,m}^{B} \cap \Omega_h \neq \emptyset. 
\]
Sets $\Gamma_h^{S,(1)}$ and $\Gamma_h^{S,(2)}$ represent, respectively, the first and second layer of ghost points.
Finally, we define $\Gamma_h^S = \Gamma_h^{S,(1)} \cup \Gamma_h^{S,(2)}$.
See Fig.~\ref{fig:stencil} (right panel) for examples of neighbour grid points 
$\mathcal{N}_{l,m}^{B}$ and $\mathcal{N}_{l,m}^{S,(k)}$ for $k=1,2$
for the ghost point $(x_l,y_m)$.

In the following, we refer to $\Gamma_h$ as either $\Gamma^S_h$ or $\Gamma^B_h$, depending on the type of internal discretizations \eqref{disc:star} or \eqref{disc:box}, respectively.
Grid points that are neither internal nor ghost points are called {\it inactive grid points}. They form the set $\mathcal{D}_h \backslash (\Omega_h \cup \Gamma_h)$. 
\subsection{Boundary points}\label{sect:bdypoint}
\begin{figure}
\begin{center}
\begin{minipage}{0.35\textwidth}
\includegraphics[width=0.99\textwidth]{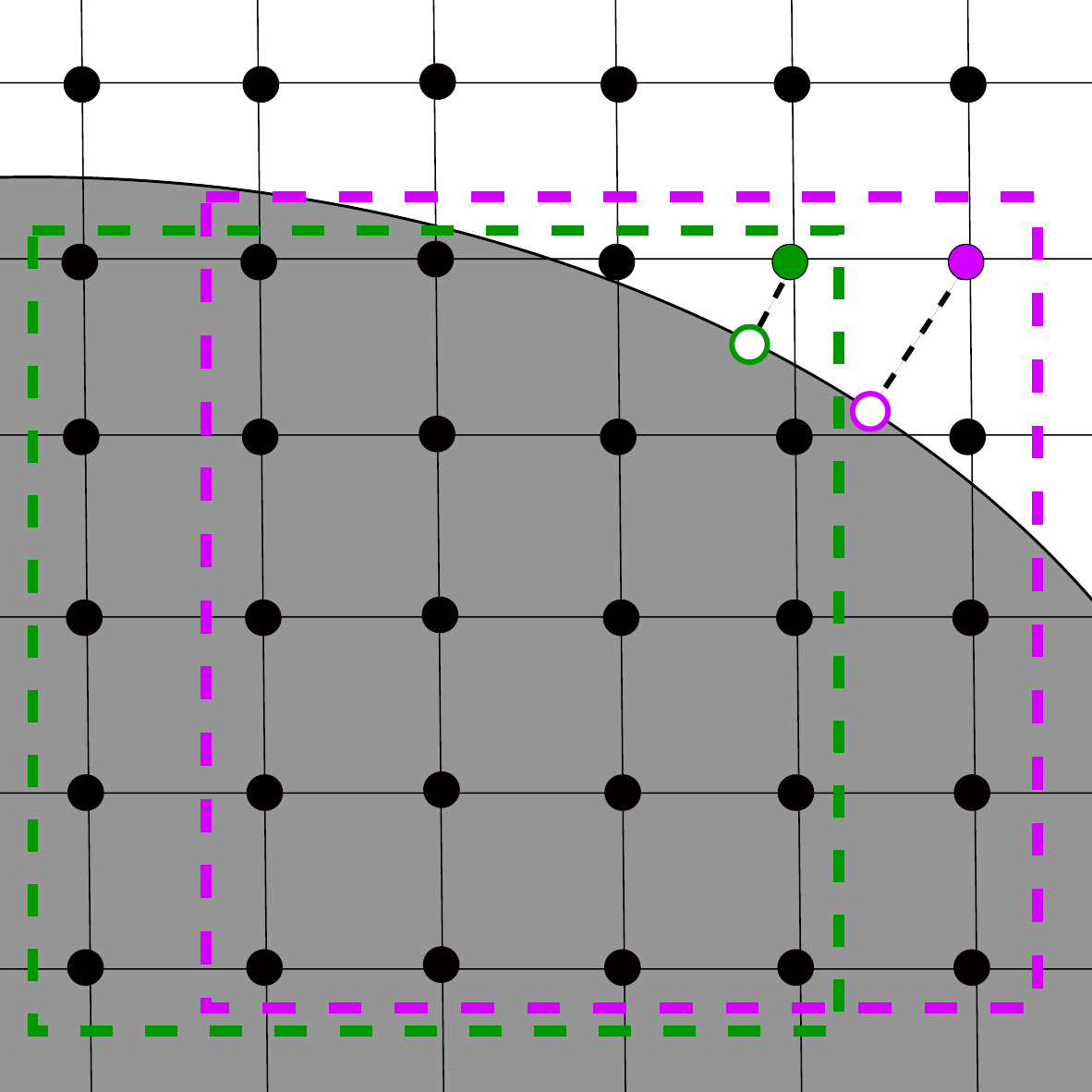}
\end{minipage}
\begin{minipage}{0.35\textwidth}
\includegraphics[width=0.99\textwidth]{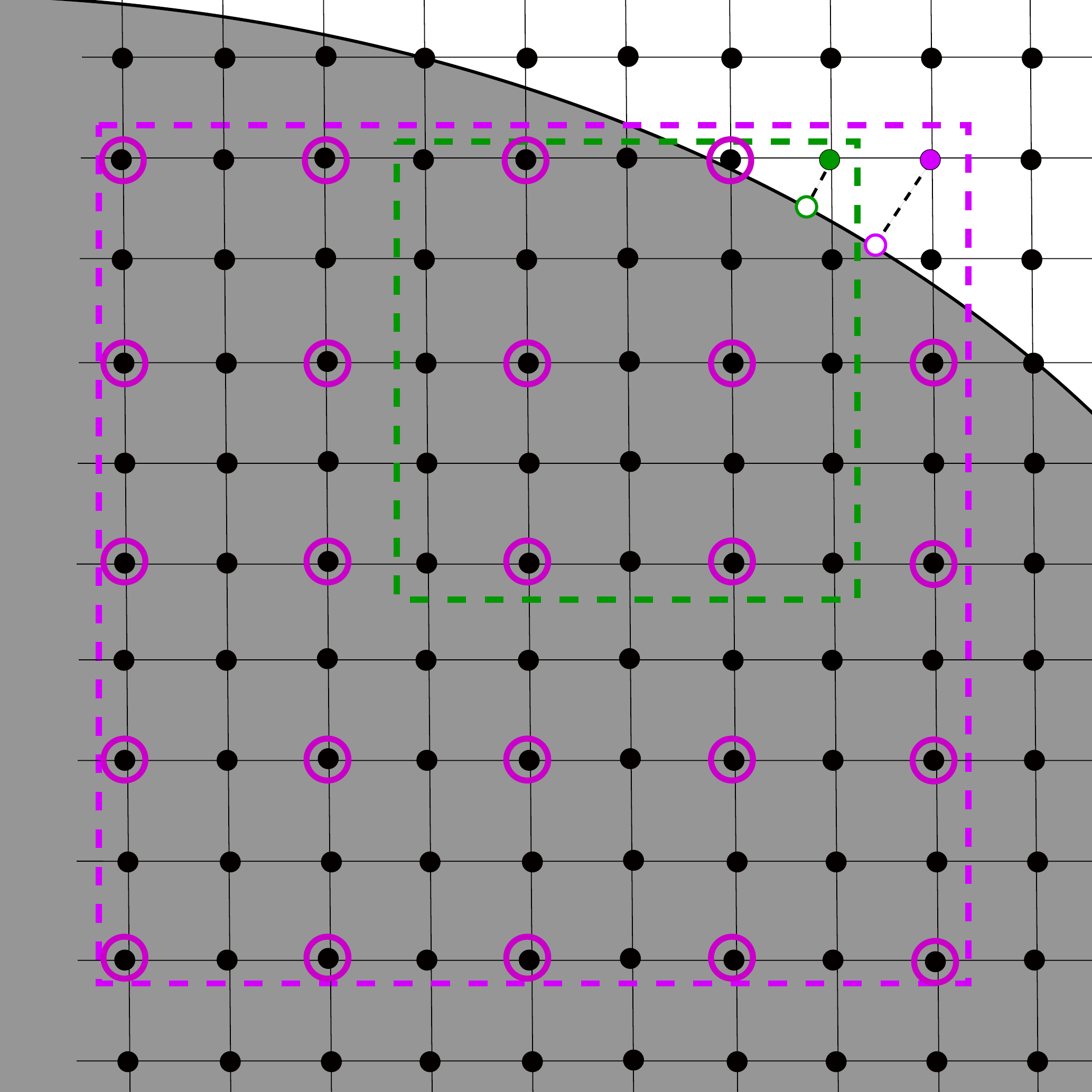}
\end{minipage}
\end{center}
\caption{
Examples of stencils $\mathcal{S}_{l,m}^{(r_x,r_y)}$ adopted in the discretization of boundary conditions to determine the ghost point values.
Left: stencils adopted in \methone and \meththree. For \methone, the green point is a ghost point of $\Gamma^{S,(1)}_h$ and the purple point is a ghost point of $\Gamma^{S,(2)}_h$. The stencils are surrounded by the dashed lines of the same colour. Both stencils have $(r_x,r_y)=(1,1)$ in Eq.~\eqref{eq:stencil}. The respective boundary points are also represented. 
Right: stencils adopted in \methtwo. The stencil for a ghost point of $\Gamma^{S,(2)}_h$ (purple point) is enlarged, meaning that $(r_x,r_y)=(2,2)$ is adopted in Eq.~\eqref{eq:stencil}. The stencil is represented by the purple point (ghost point) and the purple circles. The stencil for the green ghost point, that is a ghost point of $\Gamma^{S,(2)}_h$, is unchanged with respect to the left plot.}
\label{fig:ghostBC}
\end{figure}

For each ghost point $(x_l,y_m)$ we determine the closest boundary point $B_{l,m}$ by using the level-set function as follows \footnote{The algorithm to determine $B_{l,m}$ assumes that the level-set is known only on the grid points $\phi_{i,j}$, namely the analytical expression is not available.} 
Firstly, we divide the space into four quadrants by the Cartesian lines $x=x_l$ and $y=y_m$ and identify the quadrant $(s_x,s_y) = (\pm1,\pm1)$ in which we will search the boundary point, where $s_x=-1$ [$s_x=1$] refers to the left [right] quadrant and $s_y=-1$ [$s_y=1$] refers to the bottom [top] quadrant.
For example, $(s_x,s_y)=(1,-1)$ means that the boundary point must be found in the bottom-right quadrant. In Fig.~\ref{fig:ghostBC} we have $(s_x,s_y)=(-1,-1)$. The quadrant $(s_x,s_y)$ can be identified by looking at the opposite direction of $\nabla \phi$, the latter being orthogonal to the curve $\phi = \phi_{l,m}$. A second-order central difference discretization of $\nabla \phi$ leads to:
\[
s_x = \sgn(\phi_{l-1,m}-\phi_{l+1,m}),
\quad
s_y = \sgn(\phi_{l,m-1}-\phi_{l,m+1})
\]
where $\sgn(x)$ is the sign function defined as $\sgn(x)=1$ if $x > 0$,   $\sgn(x)=-1$ if $x < 0$ and $\sgn(0)=0$.

Then, we define the $5 \times 5-$point stencil with steps $(r_x,r_y) \in \mathbb{N}^2$ in the quadrant $(s_x,s_y)$:
\begin{equation}\label{eq:stencil}
\mathcal{S}_{l,m}^{(r_x,r_y)} =  
\left\{
(x_{l+s_x k_x r_x},y_{m+s_y k_y r_y}), \; k_x,k_y = 0,1,\ldots,4
\right\}.
\end{equation}
In this paper we will use $1 \leq r_x, r_y \leq 2$ and the specific choice of the pair $(r_x,r_y)$ will be defined later for each test.

We say that  the stencil is in the {\em upwind} direction, since the ghost point is at one corner of the stencil rather than at the centre and the rest of the stencil is taken in the opposite direction of the outer normal vector.

\subsubsection{Biquartic interpolation formulas}\label{sect:interp}
By using the $5 \times 5-$point stencil $\mathcal{S}_{l,m}^{(r_x,r_y)}$, we can interpolate a grid function $\omega_{ij}$ on an arbitrary point $(x,y)$ by evaluating the interpolant 
$\tilde{\omega}(x,y) = \sum_{p=0}^4\sum_{q=0}^4 a_{pq} x^p y^q$ where the coefficients $a_{pq}$ are determined by imposing $\tilde{\omega}(x_i,y_j) = \omega_{ij}$ for each $(x_i,y_j) \in \mathcal{S}_{l,m}^{(r_x,r_y)}$.

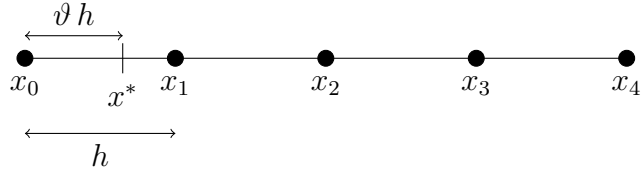
\begin{figure}
\begin{center}
\begin{tikzpicture}   
    \def\x{1.3} 
    \draw (0,0) -- (8,0);
    
    \foreach \pos/\label in {0/$x_0$,2/$x_1$,4/$x_2$,6/$x_3$,8/$x_4$} {
        \filldraw (\pos,0) circle (3pt) node[below=1mm] {\label};
    }
    
    \draw (\x, 0.2) -- (\x,-0.2) node[below] {$x^*$};

    \draw[<->] (0,0.3) -- (\x,0.3) node[midway, above] {$\theta \, h$};
    \draw[<->] (0,-1) -- (2,-1) node[midway, below] {$h$};
\end{tikzpicture}

\end{center}
\caption{1D grid nodes for quartic interpolation at $x^*$.}
\label{fig:1Dgrid}
\end{figure}

This bi-quartic interpolation can be easily computed by a tensor product of 1D quartic interpolations. In detail, the 1D quartic interpolation $\left. \tilde{\omega}(x) \right|_{x=x^*}$ of a grid function $\omega_i$ on a 5-point stencil $x_i = x_0+ih$ for $i=0,\ldots,4$ is:
\[
\left. \tilde{\omega}(x) \right|_{x=x^*} = \sum_{p=0}^4 c_p(\vartheta) \omega_{i+p}
\]
where $\vartheta = (x^*-x_0)/h$ (see Fig.~\ref{fig:1Dgrid}) and the interpolation coefficients are:

\begin{align}
c_0(\vartheta) & = \frac{1}{24}(\vartheta-4)(\vartheta-3)(\vartheta-2)(\vartheta-1), & 
c_1(\vartheta) & =  -\frac{1}{6}(\vartheta-4)(\vartheta-3)(\vartheta-2) \vartheta, \nonumber \\
c_2(\vartheta) & =  \frac{1}{4}(\vartheta-4)(\vartheta-3)(\vartheta-1) \vartheta, & 
c_3(\vartheta) & = -\frac{1}{6}(\vartheta-4)(\vartheta-2)(\vartheta-1)\vartheta, \nonumber \\
c_4(\vartheta) & =  \frac{1}{24}(\vartheta-3)(\vartheta-2)(\vartheta-1) \vartheta. \label{eq:coeffs}
\end{align}

Similarly, the first derivative is obtained by
\[
\left. \tilde{\omega}'(x) \right|_{x=x^*} = \frac{1}{h}\sum_{p=0}^4 c_p'(\vartheta) \omega_{i+p}
\]
where $\vartheta = (x^*-x_0)/h$ and the interpolation coefficients are:
\[
c_0'(\vartheta) = \frac{1}{12}(2\vartheta-5)(5 - 5 \vartheta+ \vartheta^2),\quad 
c_1'(\vartheta) = \frac{1}{6}(24 - 52\vartheta+ 27 \vartheta^2 - 4\vartheta^3), \quad
\]
\[
c_2'(\vartheta) = \frac{1}{2}(\vartheta-2)(3 - 8\vartheta+ 2\vartheta^2), \quad 
c_3'(\vartheta) = \frac{1}{6}(8 - 28\vartheta + 21\vartheta^2 - 4\vartheta^3), \quad
\]
\[
c_4'(\vartheta) = \frac{1}{12}(2\vartheta-3)(1 - 3\vartheta + \vartheta^2) 
\]

In 2D, the bi-quartic interpolation formula is
\begin{equation}\label{eq:interpD}
\tilde{\omega}(x,y) = \sum_{p,q=0}^4 c_{p}(\vartheta_x) c_{q}(\vartheta_y) \omega_{l+s_x r_x p, m+s_y r_y q}
\end{equation}
and the first partial derivatives are interpolated as
\begin{subequations} \label{derivatives}
\begin{align}
\frac{\partial \tilde{\omega}}{\partial x}(x,y) & = \frac{s_x}{h} \sum_{p,q=0}^4 c'_{p}(\vartheta_x) c_{q}(\vartheta_y) \omega_{l+s_x r_x p, m+s_y r_y q}    \\
\frac{\partial \tilde{\omega}}{\partial y}(x,y) & = \frac{s_y}{h} \sum_{p,q=0}^4 c_{p}(\vartheta_x) c'_{q}(\vartheta_y) \omega_{l+s_x r_x p, m+s_y r_y q}
\end{align}
\end{subequations}
where
\begin{equation}\label{eq:thetaext}
\vartheta_x = s_x \frac{x-x_l}{r_x\, h}, \quad \vartheta_y = s_y \frac{y-y_m}{r_y\,h}.
\end{equation}

\subsubsection{Iteration method to determine the boundary points}  
The boundary point $B_{l,m}$ is determined by the following iterative method
\[
\vec{x}^{(k+1)}=\vec{x}^{(k)} + \alpha_k \vec{p}_k, \quad \vec{x}^{(0)}=(x_l,y_m).
\]
At each step we approximate $\phi (\vec{x}^{(k)})$ and $\nabla_h {\phi}(\vec{x}^{(k)})$ by the bi-quartic interpolation $\tilde{\phi} (\vec{x}^{(k)})$ and $\nabla_h \tilde{\phi}(\vec{x}^{(k)})$ as described in Sect.~\ref{sect:interp}.
The search direction is an approximation of the unit normal vector
\begin{equation}\label{eq:searchpk}
\vec{p}^{(k)} = \frac{\nabla_h \tilde{\phi}(\vec{x}^{(k)})}{|\nabla_h \tilde{\phi}(\vec{x}^{(k)})|}
\end{equation}
and the step length is proportional to the level-set function, normalized by its gradient 
\[
\alpha^{(k)} = \varepsilon \frac{ \tilde{\phi}(\vec{x}^{(k)})}{|\nabla_h \tilde{\phi}(\vec{x}^{(k)})|}
\]
Numerically, we directly compute the product:
\[
\alpha^{(k)} \vec{p}^{(k)} = \varepsilon \, \tilde{\phi}(\vec{x}^{(k)}) \frac{\nabla_h \tilde{\phi}(\vec{x}^{(k)})}{(\nabla_h \tilde{\phi}(\vec{x}^{(k)}))^2}.
\]
The constant of proportionality must be small respect to the spatial step, say $\varepsilon = h/10$.

\begin{remark}
Selecting \(\vec{p}^{(k)} = \frac{\nabla_h \tilde{\phi}(\vec{x}^{(0)})}{|\nabla_h \tilde{\phi}(\vec{x}^{(0)})|}\) is sufficient for a signed distance function defined from the boundary to the ghost point $\vec{x}^{(0)}$. However, in practice, complications may arise. If the level-set function is not a signed distance function, if the signed distance function has not been evolved up to $\vec{x}^{(0)}$, or if the ghost point is adjacent to a concave boundary (where the gradient of the signed distance function might be discontinuous), then the chosen search direction might fail to intersect the boundary. Therefore, it is advisable to update the search direction at each iteration point as in~\ref{eq:searchpk} to accommodate these potential issues.
\end{remark}

\subsection{Enforcing boundary conditions}
For each ghost point $(x_l,y_m)$, we write a linear equation involving grid values $u_{i,j}$ on the stencil \eqref{eq:stencil}. 
The equation is obtained by enforcing the boundary condition on $B_{l,m}$ up to the desired accuracy.
Let  $\tilde{u}$ be the bi-quartic interpolant on the stencil \eqref{eq:stencil} obtained as in Sect.~\ref{sect:interp}. 
If $B_{l,m} \in \Gamma_D$, then we use the Dirichlet boundary condition and the equation is
\begin{equation}\label{eq:discD}
\tilde{u}(B_{l,m}) = g_D(B_{l,m}),
\end{equation}
otherwise (namely $B_{l,m} \in \Gamma_N$) we enforce the Neumann boundary condition:
\begin{equation}\label{eq:discN}
\nabla \tilde{u}(B_{l,m}) \cdot \frac{\nabla \phi (B_{l,m})}{|\nabla \phi (B_{l,m})|} = g_N(B_{l,m}).
\end{equation} 

\begin{remark}
The linear equation for a ghost point, either \eqref{eq:discD} or \eqref{eq:discN}, may contain other ghost values. For example, in Fig.~\ref{fig:ghostBC} (left plot), the stencil for the purple ghost point contains the green point, that is also a ghost point. Therefore, the equations for ghost points are all coupled each other. Internal and ghost equations are then solved simultaneously (see \ref{alg:LS} to assemble the linear system).
\end{remark}

\begin{remark}
We observe that Eqs.~\eqref{eq:discD} and \eqref{eq:discN} are linear equations of the grid values $u_{i,j}$ in $\mathcal{S}_{l,m}^{(r_x,r_y)} \subseteq \mathcal{D}_h$. If $\partial \Omega$ is sufficiently regular (say $\partial \Omega \in \mathcal{C}^1$) and the spatial step is sufficiently small with respect to the radius of curvature of $\partial \Omega$, then $\mathcal{S}_{l,m}^{(r_x,r_y)} \subseteq \Omega_h \cup \Gamma_h$, where $\Gamma_h$ is either $\Gamma^S_h$ or $\Gamma^B_h$ according to the chosen internal discretization \eqref{disc:star} or \eqref{disc:box}, respectively. We assume in this paper that we are always in this case.
However, if $\partial \Omega$ is only a Lipschitzian boundary, in the vicinity of kink points it may happen that $\mathcal{S}_{l,m}^{(r_x,r_y)}$ is not contained in $\Omega_h \cup \Gamma_h$ regardless of the spatial step and one has to modify the stencil accordingly, following a similar approach as in~\cite{CocoRusso:Elliptic} for the second order method.

We can label the inactive points of $\mathcal{S}_{l,m}^{(r_x,r_y)} \backslash \left( \Omega_h \cup \Gamma_h \right)$ as ghost points and add them to $\Gamma_h$. For each additional ghost point $(x_{l*},y_{m*})$, we compute again the stencil $\mathcal{S}_{l*,m*}^{(r_x,r_y)}$ and check weather it is contained in $\Omega_h \cup \Gamma_h$ or not. If not, we repeat the process again. After a few iterations, we end up with the case where the stencil $\mathcal{S}_{l,m}^{(r_x,r_y)}$ is contained in $\Omega_h \cup \Gamma_h$ for any ghost point $(x_l,y_m)$.
\end{remark}


\begin{remark}\label{remark:extrapolate_fij}
When the Mehrstellen discretization \eqref{disc:box} is adopted, the grid values $f_{ij}$ are needed on a five-point stencil centred in $(x_i,y_j)$, not only on $(x_i,y_j)$. This means that $f_{i,j}$ must be defined on some external grid points. For example, in Fig.~\ref{fig:stencil} the discretization of the Laplace equation on the grid point $(x_i,y_j)$ requires the values of $f$ on the external grid points $(x_{i+1},y_j)$ and $(x_{i},y_{j+1})$.
We do not have any information on the boundary values of $f(x,y)$, meaning that we have to extrapolate a grid function from internal to external grid points, without using any boundary condition. The extrapolation must be done to the desired accuracy, that is fourth-order accuracy in this paper. In the context of arbitrary domains represented by level-set functions, the high-order extrapolation can be obtained as in~\cite{Aslam:extrapolation}.
\end{remark}

\subsection{Linear Solver}\label{sect:LS}
The whole discretization process (at internal and ghost points) leads to a linear system in the $u_{i,j}$ for $(x_i,y_j) \in \Omega_h \cup \Gamma_h$. For simplicity, we assemble a linear system for the unknowns $u_{i,j}$ in all grid points $(x_i,y_j) \in \mathcal{D}_h$ and keep in mind that the value of $u_{i,j}$ on inactive points $\mathcal{D}_h \backslash (\Omega_h \cup \Gamma_h)$ is meaningless.
We investigate three methods in this paper:
\begin{itemize}
\item \methone - star-shaped internal discretization \eqref{disc:star} and $\Gamma_h = \Gamma^S_h$ with $(r_x,r_y)=(1,1)$ in \eqref{eq:stencil} for all ghost points
\item \methtwo - star-shaped internal discretization \eqref{disc:star} and $\Gamma_h = \Gamma^S_h$ with $(r_x,r_y)=(k,k)$ in \eqref{eq:stencil} for ghost points in $\Gamma^{S,(k)}_h$ for $k=1,2$
\item \meththree - box-shaped internal discretization (Mehrstellen) \eqref{disc:box} and $\Gamma_h = \Gamma^B_h$ with $(r_x,r_y)=(1,1)$ for all ghost points
\end{itemize}
In summary, the linear system $A \vec{u}  = \vec{b}$ is assembled by the Algorithm \ref{alg:LS}.

\begin{algorithm}
\begin{algorithmic}
\State Initialize the matrix $A$ as the identity matrix of order $(n_x+1) (n_y+1)$
 and  $\vec{b}$ as the zero column vector
\For{$(x_i,y_j) \in \Omega_h$}
\State $k = i + j n_x$
\If{ \methone  \textbf{ or} \methtwo}
\State $A(k,[k, k-1, k+1, k-n_x,k+n_x, k-2, k+2, k-2n_x, k+2n_x])$ \par 
\phantom{"  "} $= [60, -16, -16, -16, -16, 1, 1, 1, 1]/(12 h^2)$
\State $b(k) = f_{i,j}$
\ElsIf{ \meththree }
\State $A(k,[k,k-1,k+1,k-n_x,k+n_x,k-1-n_x,k-1+n_x,k+1-n_x,k+1+n_x])$ \par
\phantom{"  "} $= [20, -4, -4, -4, -4, -1, -1, -1, -1]/(6h^2)$
\State $b(k) = (8f_{i,j}+f_{i+1,j}+f_{i-1,j}+f_{i,j+1}+f_{i,j-1})/12$
\EndIf
\EndFor
\For{$(x_l,y_m) \in \Gamma_h$}
\State Find the boundary point $B_{l,m} \equiv (x_B,y_B)$ and $s_x, s_y$ by the method described in Sect.~\ref{sect:bdypoint}
\State $\vartheta_x = s_x (x_B-x_l)/(r_x \, h)$
\State $\vartheta_y = s_y (y_B-y_m)/(r_y \, h)$
\State $k = l + m\,n_x$
\If{ \methtwo \textbf{ and} $(x_l,y_m) \in \Gamma^{S,(2)}_h$}
\State $(r_x,r_y) = (2,2)$
\Else
\State $(r_x,r_y) = (1,1)$
\EndIf
\For{ $l =1, \ldots, 5$}
     \For{ $m =1, \ldots, 5$}
        \State $\texttt{stencil}(l+5(m-1))=k+s_x r_x(l-1)+s_y r_y (m-1)$
    \EndFor
\EndFor
\State $\texttt{weights}_x=[c_1(\vartheta_x), \ldots, c_5(\vartheta_x)], \quad \texttt{weights}_y=[c_1(\vartheta_y), \ldots, c_5(\vartheta_y)]$
\If{$B_{l,m} \in \Gamma_D$}
\State $b(k) = g_{D}(B_{l,m})$
\State ${\tt coeffs}={\tt weights}_{x} \cdot ({\tt weights}_{y})^{\top}$
\ElsIf{$B_{l,m} \in \Gamma_N$}
\State $b(k) = g_{N}(B_{l,m})$
\State $\texttt{weights}_{dx}=[c'_1(\vartheta_x), \ldots, c'_5(\vartheta_x)]/r_x, \quad \texttt{weights}_{dy}=[c'_1(\vartheta_y), \ldots, c'_5(\vartheta_y)]/r_y$
\State ${\tt coeffs}_{dx} = s_x {\tt weights}_y \cdot ({\tt weights}_{dx})^{\top}, \quad {\tt coeffs}_{dy} = s_y {\tt weights}_x \cdot ({\tt weights}_{dy})^{\top}$
\State $n_x= \phi({\tt stencil}) \cdot  ({\tt coeffs}_{dx})^{\top}, \quad n_y= \phi({\tt stencil}) \cdot  ({\tt coeffs}_{dy})^{\top}$
\State $n_x = n_x /\sqrt{n_x^2+n_y^2}, \quad n_y = n_y /\sqrt{n_x^2+n_y^2}$
\State ${\tt coeffs}=n_x {\tt coeffs}_{dx} +n_y {\tt coeffs}_{dy}$
\EndIf
\State $A(k,{\tt stencil}) = {\tt coeffs}$
\EndFor
\end{algorithmic}
\caption{Pseudocode to assemble the matrix $A$ and the right-hand side $b$ of the linear system (Sect.~\ref{sect:LS})}
\label{alg:LS}
\end{algorithm}

Finally, we adopt a multigrid method to solve $A \vec{u}  = \vec{b}$, similarly to the one described in~\cite{CocoRusso:Elliptic} for second order accuracy. In this paper we are interersted in investigating the accuracy of the discretization and the conditioning of the system. The efficiency of the solver is not within the scope of this study, although we experienced a similar performance to what is observed for the second order accurate method in~\cite{CocoRusso:Elliptic, coco2018second}.

\section{Numerical results}\label{sect:numtests}
\begin{figure}
\begin{center}
\begin{minipage}{0.4\textwidth}
\includegraphics[width=0.99\textwidth]{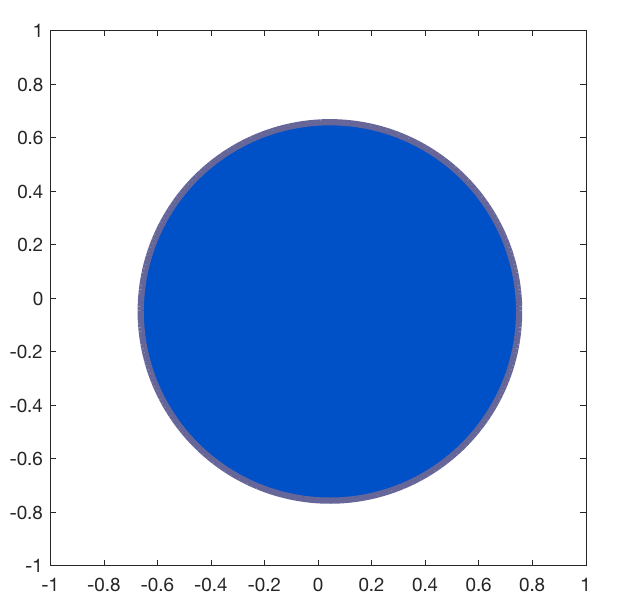}
\end{minipage}
\begin{minipage}{0.4\textwidth}
\includegraphics[width=0.99\textwidth]{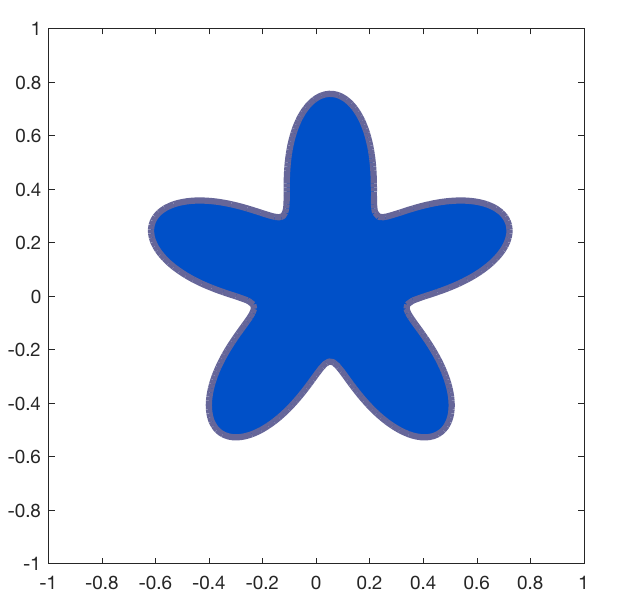}
\end{minipage}
\end{center}
\caption{
Domains used in numerical tests (Sect.~\ref{sect:numtests}).
}
\label{fig:domains}
\end{figure}

We perform numerical tests on two different domains (see Fig.~\ref{fig:domains})): a circle and a flower-shaped domain.
The level-set representations $\phi_C$ and $\phi_F$ of the two domains are:
\begin{multline}\label{dom2D:circle}
\text{\sc{(Circle)}} 
\quad
 \phi_C(x,y) = \sqrt{(x-x_0)^2 + (y-y_0)^2} - R
\text{ with } \; x_0=\sqrt{2}/10, \; y_0=-\sqrt{3}/20, \; R=\sqrt{5}/3
\end{multline}
\begin{multline}\label{dom2D:flower}
\text{\sc{(Flower)}} 
\quad
\phi_F(x,y) = r-r_1
     -\frac{r_2}{r^5} \left( (y-y_0)^5+5(x-x_0)^4 (y-y_0)-10(x-x_0)^2 (y-y_0)^3 \right), \\
    \text{ with } r = \sqrt{x^2+y^2}, \quad r_1=0.5, \quad r_2=0.2.
\end{multline}

The flower-shaped domain can also be described by the parametric representation:
\begin{equation*}
\left\{
\begin{aligned}
x(\theta) &= r(\gamma) \cos(\gamma) + x_0 \\
y(\theta) &= r(\gamma) \sin(\gamma) + y_0
\end{aligned}
\right.
\end{equation*}
with $\gamma \in [0,2\pi]$ and $r(\gamma) = r_1 + r_2 \sin(5 \gamma)$.

For each test we choose an exact solution $u_\text{EXA}(x,y)$ and then compute analitically $f$, $g_D$ and $g_N$ from \eqref{eq:modprob}. We observe that the analytical expression of $f(x,y)$ can be used to define grid values outside the domain. However, we assume that $f_{ij}$ is only known on internal grid points and then we extrapolate it outside $\Omega$ according to Remark~\ref{remark:extrapolate_fij}.

The partition of the boundary into $\Gamma_D$ and $\Gamma_N$ is performed in the following way in all numerical tests: a boundary point of $\partial \Omega$ is in $\Gamma_D$ if $x \leq 0$, otherwise it is in $\Gamma_N$.

The relative numerical error is computed as
\[
e_h = \frac{\left\| u_h - u_\text{EXA} \right\|_p}{\left\| u_\text{EXA} \right\|_p}
= \left( \frac{\sum_{(x_i,y_j) \in \Omega_h \cup \Gamma_h} \left| u_{i,j} - u_\text{EXA}(x_i,y_j) \right|^p }{\sum_{(x_i,y_j) \in \Omega_h \cup \Gamma_h} \left| u_\text{EXA}(x_i,y_j) \right|^p } \right)^{1/p},
\quad \text{ with } p=1, \infty,
\]
where $u_h$ is the numerical solution defined on the grid nodes of $\Omega_h \cup \Gamma_h$. 
We also approximate the gradient of the solution by a fourth-order accurate finite-difference formula:
\begin{equation}\label{eq:discnabla}
\nabla u_h = \left( \frac{\partial u_h}{\partial x}, \frac{\partial u_h}{\partial y} \right)=
\left(\frac{-u_{i+2,j}+8u_{i+1,j}-8u_{i-1,j}+u_{i-2,j}}{12 h},
\frac{-u_{i,j+2}+8u_{i,j+1}-8u_{i,j-1}+u_{i,j-2}}{12 h} \right) + O(h^4).
\end{equation}
Then, we compute the error on the gradient of the solution by:
\[
e^{\nabla u}_h = \frac{\left\| \nabla u_h - \nabla u_\text{EXA} \right\|_p}{\left\| \nabla u_\text{EXA} \right\|_p}
= \left( \frac{\sum_{(x_i,y_j) \in \Omega_h} \left| \nabla u_{i,j} - \nabla u_\text{EXA}(x_i,y_j) \right|^p }{\sum_{(x_i,y_j) \in \Omega_h} \left| \nabla u_\text{EXA}(x_i,y_j) \right|^p } \right)^{1/p},
\quad \text{ with } p=1, \infty,
\]
where $\left| (v_1,v_2) \right| = \sqrt{v_1^2+v_2^2}$ is the modulus of a two-dimensional vector.

\subsection{\methone}\label{sect:numtest_meth1}
\begin{figure}
\begin{center}
\begin{minipage}{0.4\textwidth}
\includegraphics[width=0.99\textwidth]{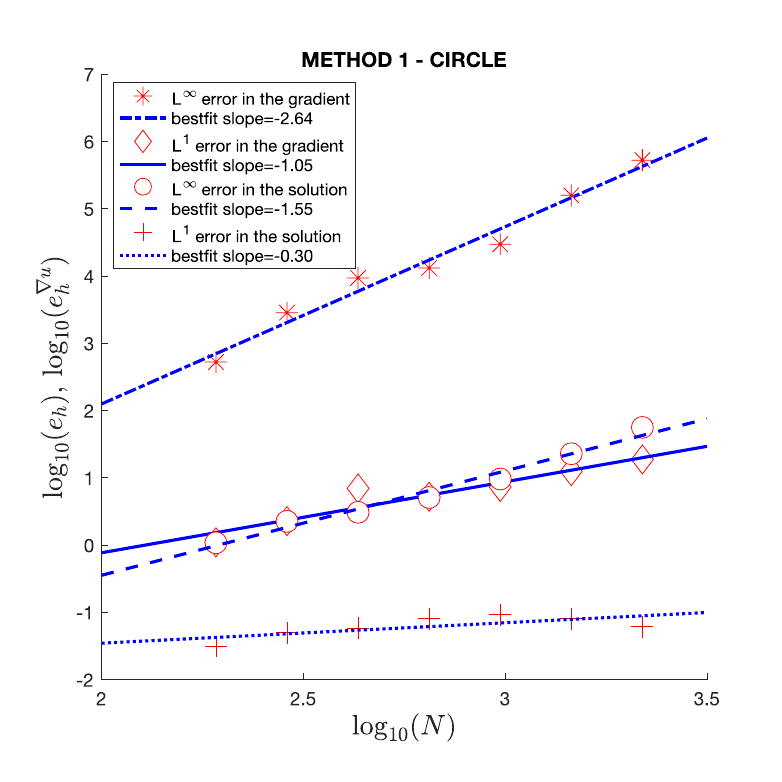}
\end{minipage}
\begin{minipage}{0.4\textwidth}
\includegraphics[width=0.99\textwidth]{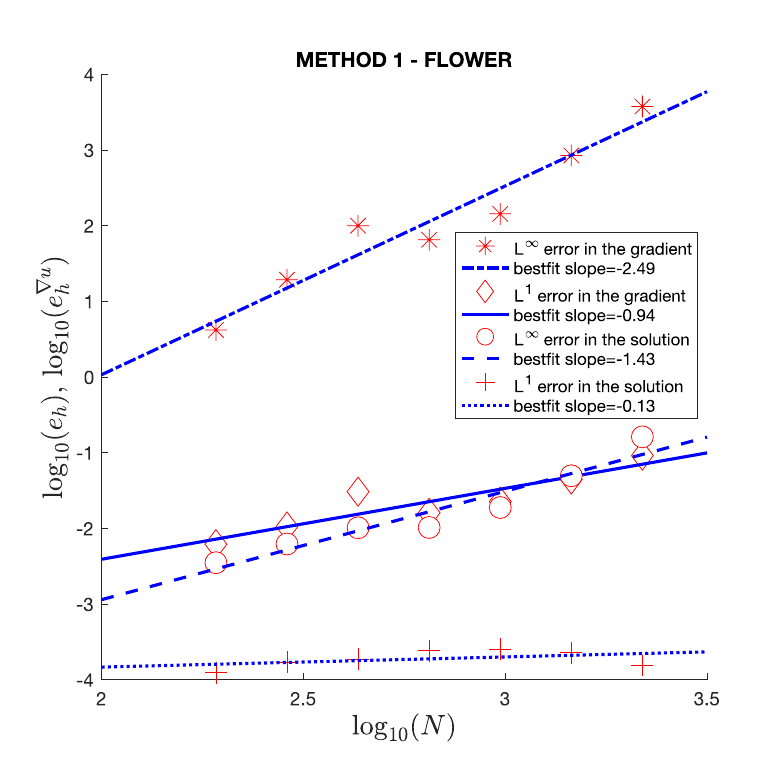}
\end{minipage}
\end{center}
\caption{
\methone: numerical error of the solution and the gradient in $L^1$ and $L^\infty$ norms vs different values of $N$, for both the circular (left) and the flower-shaped (right) domains. The slope of the bestfit lines provide an approximation of the accuracy orders, shown in the legend.
}
\label{fig:acc_METH1}
\end{figure}

\begin{figure}
\begin{center}
\begin{minipage}{0.40\textwidth}
\includegraphics[width=0.99\textwidth]{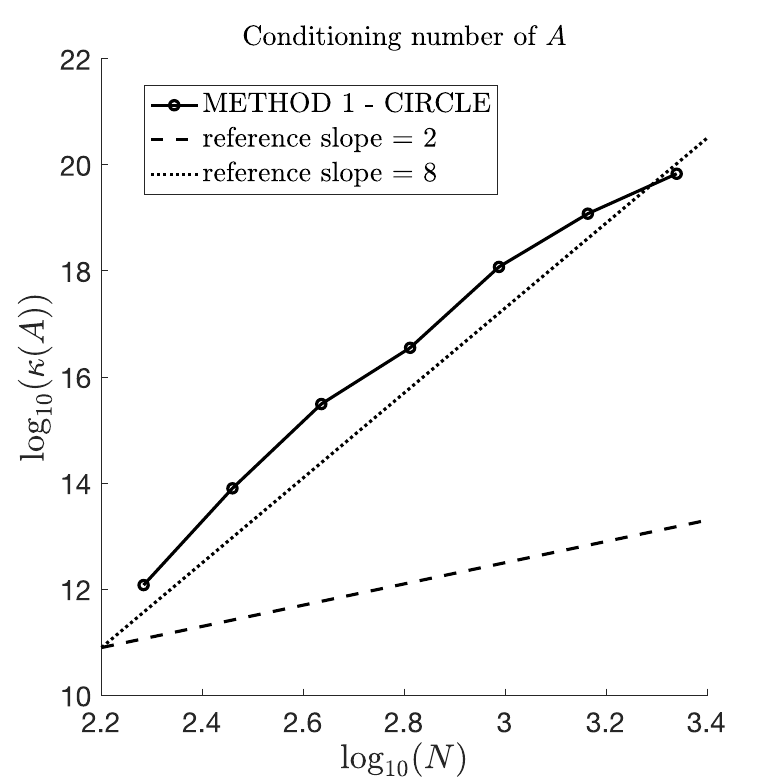}
\end{minipage}
\begin{minipage}{0.40\textwidth}
\includegraphics[width=0.99\textwidth]{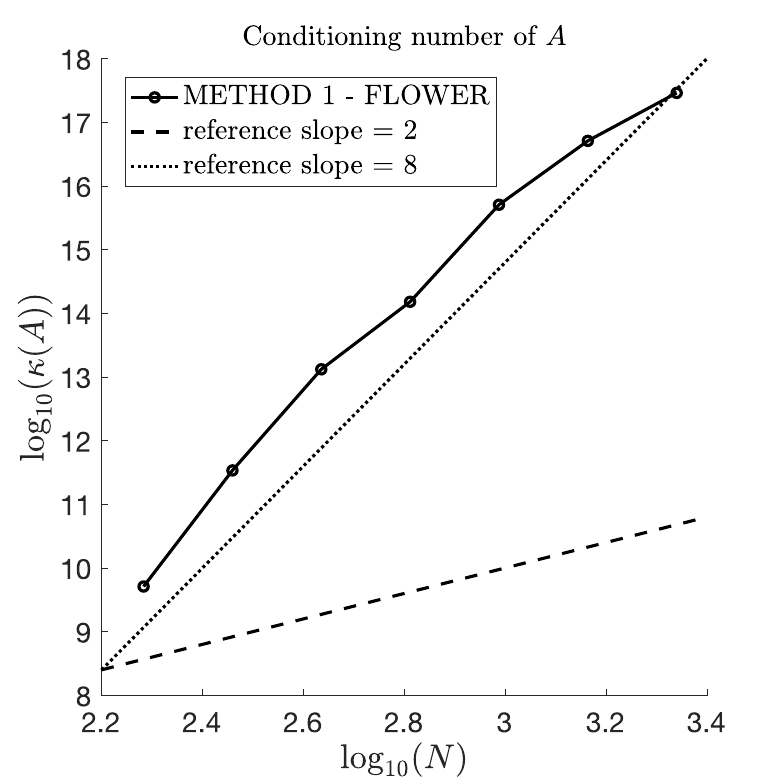}
\end{minipage}
\end{center}
\caption{
\methone: condition number of the matrix $A$ for different values of $N$ for the circular (left) and the flower-shaped (right) domains. Reference slopes of $O(N^2)$ and $O(N^8)$ are represented by dashed and dotted lines, respectively.}
\label{fig:cond_METH1}
\end{figure}

We choose $u_\text{EXA}=\sin(2x)\sin(5y)$ for the circular domain and $u_\text{EXA}=\log(1+3xy)$ for the flower-shaped domain. The errors are plotted in Fig.~\ref{fig:acc_METH1}, showing that the method is not consistent, namely the numerical error does not even converge to the exact solution when $N$ goes to infinity. Indeed, it is an increasing function of $N$. This behavior is explained by the condition number of the matrix $A$, that is very large (close to $10^{16}$ and then almost no digits are accurate in the numerical solution), resulting in an ill-conditioned linear system, as shown in Fig.~\ref{fig:cond_METH1}. 
This is caused by the presence of possible ill-conditioned interpolations. For example, in \eqref{eq:interpD} for Dirichlet boundary conditions \eqref{eq:discD} related to the second layer of ghost-points $\Gamma^{S,(2)}_h$,
if we solve \eqref{eq:discD} for the ghost value $u_{l,m}$, we obtain
\[
u_{l,m} = \frac{1}{c_{0}(\theta_x)c_{0}(\theta_y)}\left( g_D(B_{l,m}) - \sum_{\substack{(p,q)\in \left\{0,\ldots,4\right\}^2 \\ (p,q)\neq(0,0)}}  c_{p}(\vartheta_x) c_{q}(\vartheta_y) u_{l+s_x r_x p, m+s_y r_y q}\right).
\]
This is ill-conditioned when $c_{0}(\theta_x)c_{0}(\theta_y) \approx 0$, namely when $\theta_x$ or $\theta_y$  are close to 1 (see \eqref{eq:coeffs}). This is observed, for example, in Fig.~\ref{fig:ghostBC} for the purple ghost point (that belongs to $\Gamma^{S,(2)}_h$), since  its boundary projection is close to the grid point $(x_{l-1},y_{m-1})$. To overcome this issue, \methtwo enlarges the stencil for ghost points of $\Gamma^{S,(2)}_h$ by using $(r_x,r_y) = (2,2)$ in \eqref{eq:stencil}. In this way, $\theta_x$ and $\theta_y$ are smaller than $1/2$ (from \eqref{eq:thetaext}) and the interpolation is always well-conditioned.

The ill-conditioned behavior can also be explained in terms of small diagonal terms of the matrix $A$ (compared to off-diagonal terms of the same row). In fact, in the case of Dirichlet boundary condition, the diagonal term is $c_{0}(\theta_x)c_{0}(\theta_y)$ and the off-diagonal terms are $c_{p}(\theta_x)c_{q}(\theta_y)$ for $(p,q)\in \left\{0,\ldots,4\right\}^2$, $(p,q)\neq(0,0)$. If $\theta_x$ or $\theta_y$ are close to 1, then the diagonal term $c_{0}(\theta_x)c_{0}(\theta_y) \approx 0$. 

We note that the issue is observed in several methods where the domain is embedded in a fixed grid that does not conform the boundary, falling under the category of unfitted boundary methods, such as other ghost point methods~\cite{gibou2002second} or cut-cell methods (where the issue is caused by the presence of small cut cells~\cite{schneiders2013accurate}).

\subsection{\methtwo}\label{sect:numtest_meth2}
\begin{figure}
\begin{center}
\begin{minipage}{0.40\textwidth}
\includegraphics[width=0.99\textwidth]{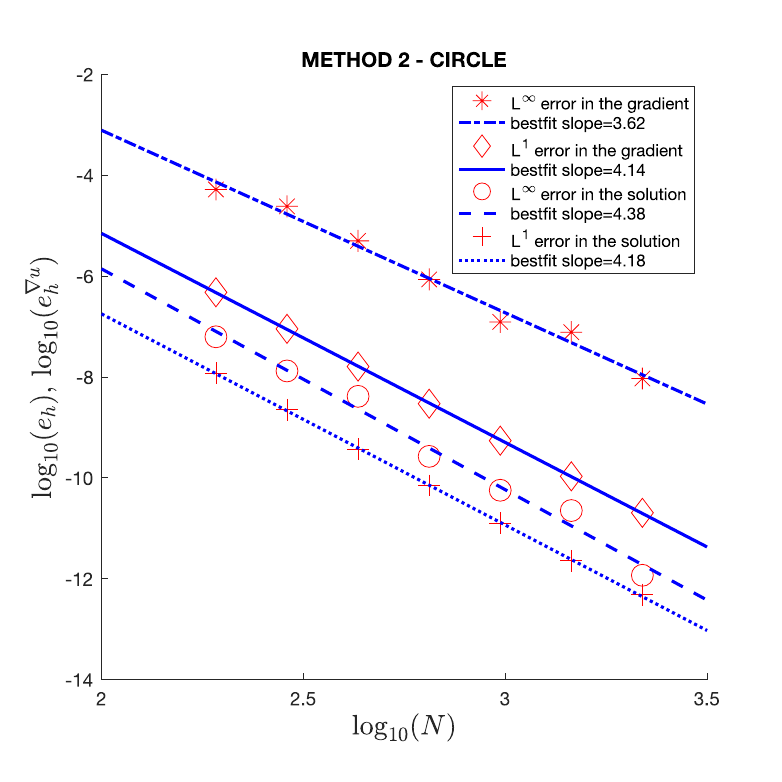}
\end{minipage}
\begin{minipage}{0.40\textwidth}
\includegraphics[width=0.99\textwidth]{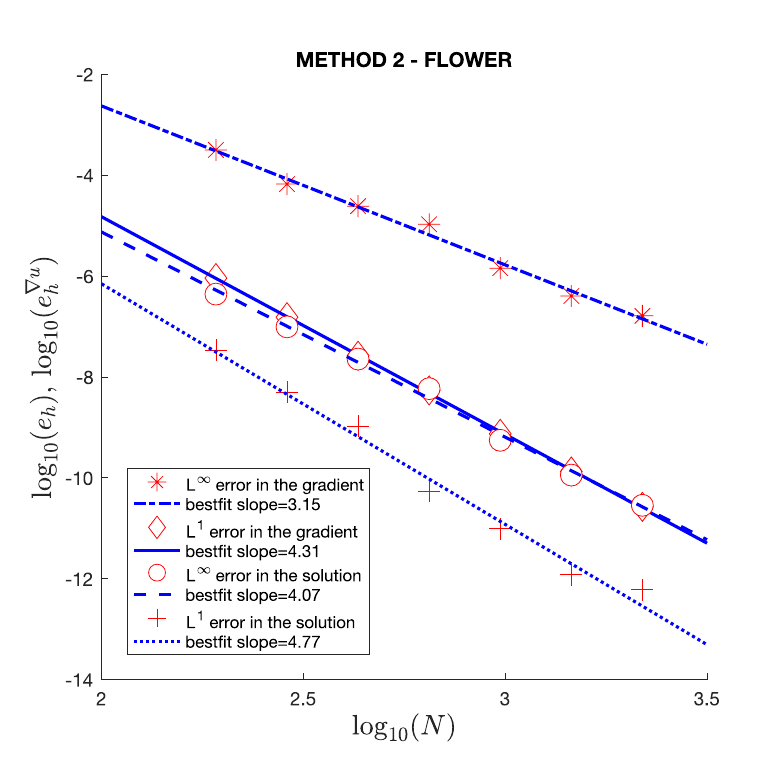}
\end{minipage}
\end{center}
\caption{
\methtwo: numerical error of the solution and the gradient in $L^1$ and $L^\infty$ norms vs different values of $N$, for both the circular (left) and the flower-shaped (right) domains. The slope of the bestfit lines provide an approximation of the accuracy orders, shown in the legend.
}
\label{fig:acc_METH2}
\end{figure}

We perform the same tests as in~\ref{sect:numtest_meth1}, except that we use en enlarged stencil for the ghost points of the second layer $\Gamma^{S,(2)}_h$, namely $(r_x,r_y)=(2,2)$ in \eqref{eq:stencil}, so that $\vartheta_x$ and $\vartheta_y$ are always smaller than $1/2$ and then $c_0(\vartheta_x)c_0(\vartheta_y)$ is not close to zero. We observe that the method is fourth-order accurate in the solution and its gradient (Fig.~\ref{fig:acc_METH2}) and the condition number is improved with respect to~\ref{sect:numtest_meth1}, scaling as $\kappa(A) = O(N^2)$ (Fig.~\ref{fig:cond_METH2and3}), that is the expected scaling of the discrete Laplacian matrix, meaning that the ghost point approach does not affect the overall stability of the method.

\subsection{\meththree}\label{sect:numtest_meth3}
\begin{figure}
\begin{center}
\begin{minipage}{0.40\textwidth}
\includegraphics[width=0.99\textwidth]{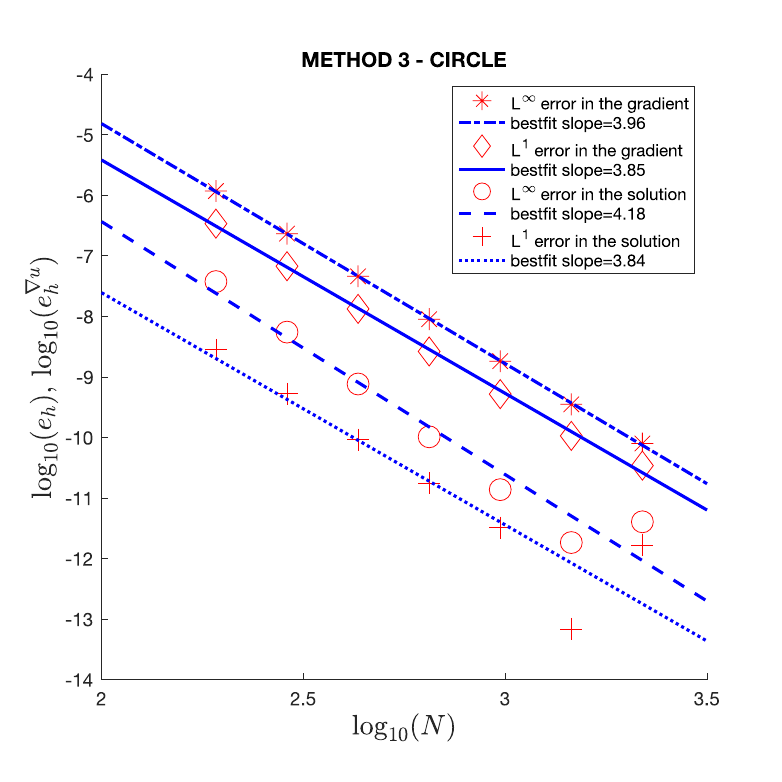}
\end{minipage}
\begin{minipage}{0.40\textwidth}
\includegraphics[width=0.99\textwidth]{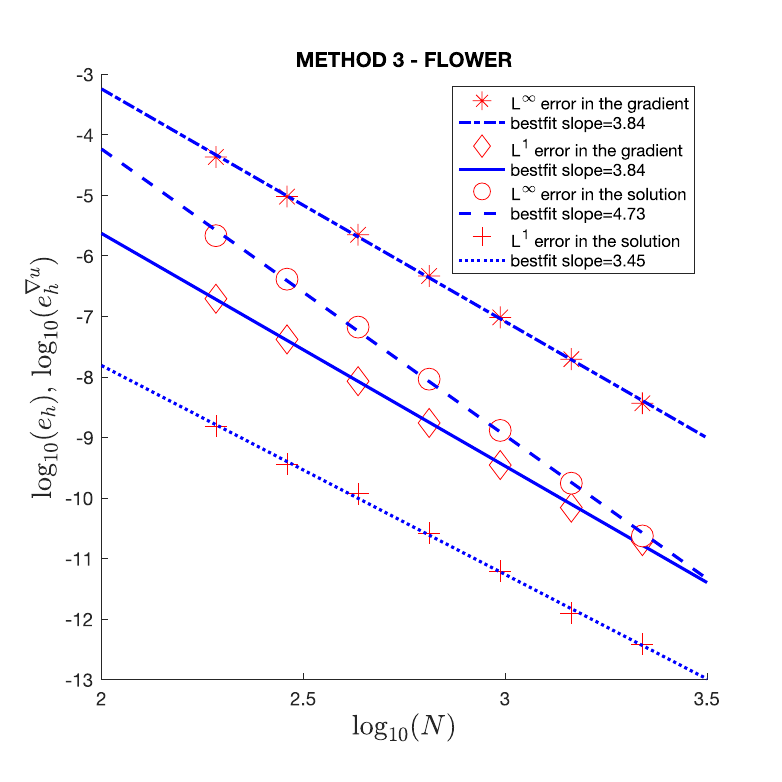}
\end{minipage}
\end{center}
\caption{
\meththree: numerical error of the solution and the gradient in $L^1$ and $L^\infty$ norms vs different values of $N$, for both the circular (left) and the flower-shaped (right) domains. The slope of the bestfit lines provide an approximation of the accuracy orders, shown in the legend.
}
\label{fig:acc_METH3}
\end{figure}

\begin{figure}
\begin{center}
\begin{minipage}{0.40\textwidth}
\includegraphics[width=0.99\textwidth]{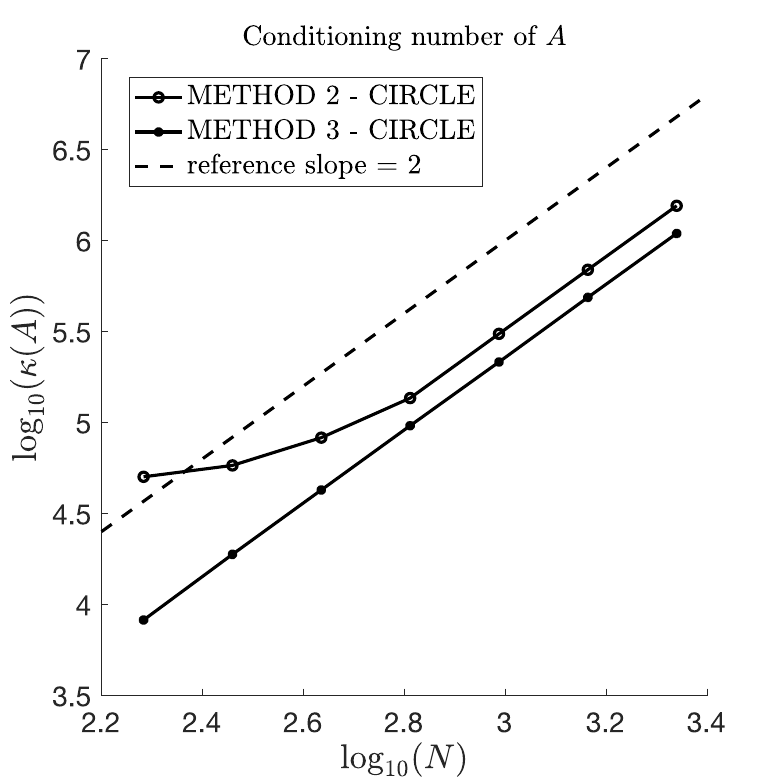}
\end{minipage}
\begin{minipage}{0.40\textwidth}
\includegraphics[width=0.99\textwidth]{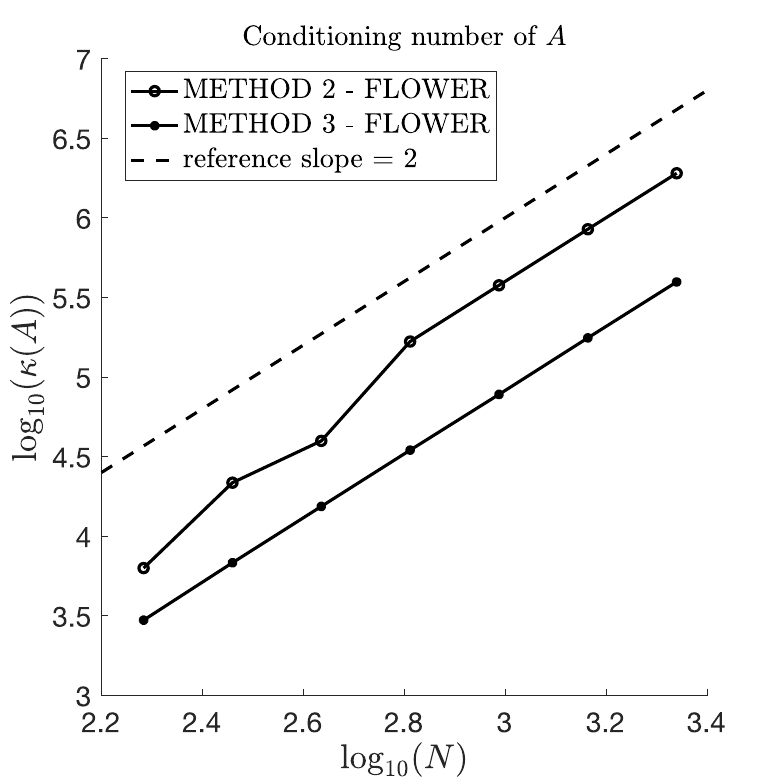}
\end{minipage}
\end{center}
\caption{
\methtwo and \meththree: condition number of the matrix $A$ for different values of $N$ for the circular (left) and the flower-shaped (right) domains. A reference slope of $O(N^2)$ is represented by a dashed line.
}
\label{fig:cond_METH2and3}
\end{figure}

We implement the box-shaped internal discretization (Mehrstellen) and plot the errors in Fig.~\ref{fig:acc_METH3}.
We observe that the central difference formula based on a stencil centered in $(x_i,y_j)$ to approximate the gradient \eqref{eq:discnabla} cannot be computed for those internal grid points that are close to the boundary, since the second layer of external points is made by inactive points in this method. Therefore, for those points we use a different fourth-order formula by selecting a different stencil, maintaining the same accuracy. For example, the discretization of ${\partial u_h}/{\partial x}$ on the grid point $(x_i,y_j)$ of Fig.~\ref{fig:stencil} is:
\[
 \frac{\partial u_h}{\partial x} =
\frac{-u_{i-3,j}+6u_{i-2,j}-18u_{i-1,j}+10u_{i,j}+3u_{i+1,j}}{12 h}+ O(h^4).
\]
instead of the one used in~\eqref{eq:discnabla}.
Like \methtwo~of Sect.~\ref{sect:numtest_meth2}, \meththree is fourth-order accurate in the solution and its gradient. The linear system is well conditioned and $\kappa(A) = O(N^2)$ as shown in Fig.~\ref{fig:cond_METH2and3}.

\subsection{Comparison of the methods}
\methone is not consistent due to the ill-conditioned matrix, and therefore is not further discussed in this section. \methtwo and \meththree are both fourth-order accurate in the solution and its gradient, say that the errors are asymptotically scaling as $C_2 h^4$ and $C_3 h^4$ for \methtwo and \meththree respectively.
However, \meththree is considerably more accurate than \methtwo meaning that $C_3 < C_2$.
Condition number of the matrix $A$ is very similar for both methods and scales as expected, $\kappa(A) = O(N^2)$, although \meththree is slightly more stable.

The improved accuracy and condition of \meththree are justified by the adoption of a better stencil for the discretization of the Laplace operator, since the box-stencil is more compact than the star-stencil. This is also preferred by the ghost-point approach, since ghost values are only required in the vicinity of the boundary, while the star-shaped stencil needs the definition of grid values on a farther second layer of ghost points.

The main drawback of \meththree lies in the requirement that the source term $f$ must be evaluated  on a $5$-point star-stencil (Eq.~\eqref{disc:box}), as highlighted in Remark~\ref{remark:extrapolate_fij}.
For internal grid points that are close to the boundary, this means that $f$ must be evaluated also on external grid points, therefore a suitable technique must be adopted to extrapolate $f$ from internal to external points. The ghost-point approach adopted in this paper to define the ghost values of the solution outside the domain cannot be employed to extrapolate $f$, since no boundary conditions are available for the source term. In this paper, the artificial solution $u_\text{EXA}$ automatically defines an analytical expression for the source term $f$ that can be used to evaluate the source terms in points outside of the domain. In real applications this might not possible, since usually only internal grid values are available. 
Therefore, one has to adopt some other techniques to extrapolate $f$ from internal to external grid points to the desired accuracy order. To this purpose, a very efficient approach is the one proposed by Aslam in \cite{Aslam:extrapolation}, where grid values can be extrapolated outside of a domain with arbitrary geometry (defined by a level-set function) up to the desired accuracy.
\section{Conclusion}
In this paper fourth order extensions of the finite difference ghost-point method introduced in~\cite{CocoRusso:Elliptic} have been presented. 
The methods are based on two main ingredients: 
\begin{itemize}
\item finite difference discretization on the internal nodes;
\item interpolation based on an upwind stencil to apply the boundary condition at the orthogonal projection of each ghost point onto the boundary. 
\end{itemize}
To attain fourth-order accuracy not only for the function $u$ but also for its gradient, it is necessary to employ biquartic interpolation. This is because bicubic interpolation only ensures fourth-order accuracy for the function and third-order accuracy for the gradient.

Three methods, which are theoretically fourth-order accurate for both the function and its gradient, are considered. However, one of these methods exhibits such poor conditioning that it fails to converge. The remaining two methods demonstrate fourth-order accuracy for both the function and the gradient. Among these, the method utilizing Mehrstellen discretization achieves superior accuracy and is more practical, as it requires only a single layer of ghost points.

This paper primarily examines the accuracy of the method while setting aside several critical aspects, such as efficiency and a detailed convergence analysis of the schemes.
Concerning efficiency, a key advantage highlighted in \cite{CocoRusso:Elliptic} is the compatibility of the method with a geometric multigrid approach for solving the linear system, which significantly enhances its efficiency. In \cite{coco2023ghost}
the efficiency of the multigrid method has been improved by optimizing the relaxation parameters used during ghost value iteration. This optimization specifically aims to smooth the residual in the tangential direction along the boundary.

The convergence analysis for rectangular domains whose boundaries do not align with Cartesian grid lines has been investigated in \cite{coco2023spectral}.

We expect that similar strategies regarding efficiency and convergence analysis can be effectively applied to the fourth-order ghost method discussed in this paper. Efforts to pursue these enhancements are actively underway.

\section*{Acknowledgement}
This work has been supported by the Spoke 1 Future HPC \& Big Data of the Italian Research Center on High-Performance Computing, Big Data and Quantum Computing (ICSC) funded by MUR Missione 4 Componente 2 Investimento 1.4: Potenziamento strutture di ricerca e creazione di “campioni nazionali di R \&S (M4C2-19)” - Next Generation EU (NGEU).

The work of A.C.~has been supported from Italian Ministerial grant PRIN 2022 “Efficient numerical schemes and optimal control methods for time-dependent partial differential equations”, No. 2022N9BM3N - Finanziato dall’Unione europea - Next Generation EU – CUP: E53D23005830006.

The work of A.C.~and G.R.~has been supported from Italian Ministerial grant PRIN 2022 PNRR “FIN4GEO: Forward and Inverse Numerical Modeling of hydrothermal systems in volcanic regions with application to geothermal energy exploitation”, No. P2022BNB97 - Finanziato dall’Unione europea - Next Generation EU – CUP: E53D23017960001.

A.C.~and G.R.~are members of the Gruppo Nazionale Calcolo Scientifico-Istituto Nazionale di Alta Matematica (GNCS-INdAM).

We acknowledge the CINECA award under the ISCRA initiative (ISCRA C project BCMG, code HP10C7YOPZ), for the availability of high-performance computing resources and support.

\addcontentsline{toc}{chapter}{References}
\bibliographystyle{abbrv}
\bibliography{bibliography}
\end{document}